\documentclass[a4paper,12pt]{amsart}\usepackage{amsfonts,amssymb,oldgerm,enumerate,mathrsfs,bbm,hyperref,rotating}

\UseRawInputEncoding

\def\C{\Bbb{C}}

\def\P{\mathbb{P}}\def\R{\Bbb{R}}\def\Z{\mathbb{Z}}

\def\di{\partial}
\def\bl{\langle}\def\br{\rangle}
\def\oplusl{\mathop\oplus\limits}

\def\supl{\sup\limits}\def\infl{\inf\limits}

\newcommand{\quots}[2]{{\footnotesize\left.\raisebox{0.4ex}{$#1$}\! / \!\raisebox{-0.4ex}{$#2$}\right.}}
\newcommand{\quot}[2]{{\left.\raisebox{1.6ex}{\footnotesize$#1$}  \!\!\!\!{\ensuremath\diagup}\!\!\raisebox{-1ex}{\footnotesize$#2$}\right.}}
\renewcommand{\stackrel}[2]{\ \lower 0.2ex \hbox{$\mathrel{\mathop{#2}\limits^{#1}}$}\ }

\def\tDe{{\tilde{\De}}}\def\teta{\tilde{\eta}}
\def\tPhi{\tilde{\Phi}}
 \def\tp{\tilde{p}}\def\tpi{{\tilde\pi}}

\def\cU{\mathcal U}

\def\tx{\tilde{x}}\def\tX{{\tilde{X}}}\def\tY{{\tilde{Y}}}

\def\hx{\hat{x}}\def\hz{\hat{z}}

\def\ga{\gamma}\def\be{\beta}\def\De{\Delta}
\def\ep{\epsilon}\def\om{\omega}\def\Om{\Omega}
\def\si{\sigma}

\def\cO{\mathcal{O}}\def\cS{\mathcal{S}}

\def\uf{{\underline{f}}}\def\uh{{\underline{h}}}

\def\ux{{\underline{x}}}

\def\empty{\varnothing}

\newcommand{\bbm}{\begin{bmatrix}}\newcommand{\ebm}{\end{bmatrix}}
\newcommand{\ber}{\begin{array}{l}}\newcommand{\eer}{\end{array}}
\newcommand{\bpm}{\begin{pmatrix}}\newcommand{\epm}{\end{pmatrix}}
\newcommand{\bM}{\begin{matrix}}\newcommand{\eM}{\end{matrix}}
\newcommand{\bee}{\begin{enumerate}}\newcommand{\eee}{\end{enumerate}}
\newcommand{\bei}{\begin{itemize}}\newcommand{\eei}{\end{itemize}}

\def\wrt{with respect to }\def\sset{\subset}\def\sseteq{\subseteq}\def\smin{\setminus}

\newcommand{\beq}{\begin{equation}}\newcommand{\eeq}{\end{equation}}
\newtheorem{Lemma}{Lemma}[section]\newcommand{\bel}{\begin{Lemma}}\newcommand{\eel}{\end{Lemma}}
\newtheorem{Example}[Lemma]{Example}\newcommand{\bex}{\begin{Example}\rm}\newcommand{\eex}{\end{Example}}
\newtheorem{Proposition}[Lemma]{Proposition}\newcommand{\bprop}{\begin{Proposition}}\newcommand{\eprop}{\end{Proposition}}
\newtheorem{Property}[Lemma]{Property}\newcommand{\bproperty}{\begin{Property}}\newcommand{\eproperty}{\end{Property}}
\newtheorem{Definition-Proposition}[Lemma]{Definition-Proposition}

\def\bpr{~\\{\em Proof.\ }}
\newcommand{\epr}{{\hfill\ensuremath\blacksquare}\\}
\newtheorem{Theorem}[Lemma]{Theorem}\newcommand{\bthe}{\begin{Theorem}}\newcommand{\ethe}{\end{Theorem}}
\newtheorem{Definition}[Lemma]{Definition}\newcommand{\bed}{\begin{Definition}}\newcommand{\eed}{\end{Definition}}
\newtheorem{Remark}[Lemma]{Remark}\newcommand{\beR}{\begin{Remark}\rm}\newcommand{\eeR}{\end{Remark}}
\newtheorem{Corollary}[Lemma]{Corollary}\newcommand{\bcor}{\begin{Corollary}}\newcommand{\ecor}{\end{Corollary}}
\newcommand{\bet}{\begin{tabular}{cccccccc}}\newcommand{\eet}{\end{tabular}}

\newcommand{\isom}{\xrightarrow[\,\smash{\raisebox{1.15ex}{\ensuremath{\scriptstyle\sim}}}\,]{}}

\title[]{D\MakeLowercase{etecting fast vanishing loops in complex-analytic germs}\\
 \MakeLowercase{(and detecting germs that are inner metrically conical)}}
\author[]{D\MakeLowercase{mitry}  K\MakeLowercase{erner and} R\MakeLowercase{odrigo} M\MakeLowercase{endes}}

\thanks{We were supported by the Israel Science Foundation,  grants No.  1910/18 and 1405/22}

\address{Dmitry Kerner: Department of Mathematics, Ben Gurion University of the Negev, P.O.B. 653, Be'er Sheva 84105, Israel.}
\email{dmitry.kerner@gmail.com}

\address{Rodrigo Mendes: Instituto de ci\^encias exatas e da natureza, Universidade de Integra\c{c}\~ao Internacional da Lusofonia Afro-Brasileira (unilab), Campus dos Palmares, Cep. 62785-000. Acarape-Ce, Brasil and  Departament of Mathematics, Ben Gurion University of the Negev, P.O.B. 653, Be'er Sheva 84105, Israel.}
\email{rodrigomendes@unilab.edu.br}


\date{\today\ \  filename: \jobname.tex}

\begin{document}\setcounter{secnumdepth}{6} \setcounter{tocdepth}{1}

 \begin{abstract}
  Let $X$ be a reduced complex-analytic germ of pure dimension $n\ge2,$ with arbitrary singularities (not necessarily normal or complete intersection). Various homology cycles on $Link_\ep[X]$ vanish at different speeds when $\ep\to0.$
  We give a  condition ensuring fast vanishing loops on $X.$  The condition is in terms of the discriminant
  and  the covering data for ``convenient" coverings $X\to (\C^n,o)$.   No resolution of singularities is involved.

    For surface germs ($n=2$) this condition becomes necessary and sufficient.

A corollary for surface germs that are strictly complete intersections detects fast loops via singularities of the projectivized tangent cone of $X.$

Fast loops are the simplest obstructions for $X$ to be inner metrically conical. Hence we get simple necessary conditions to the IMC property.
For normal surface germs these conditions are also sufficient.

We give numerous classes of non-IMC germs  and    IMC germs.
 \end{abstract}
 \maketitle
\tableofcontents

\section{Introduction}
\subsection{} Take a complex-analytic germ $X \sset (\C^N,o).$
 The first step in its visualization is Milnor's Conic structure theorem:
 {\em  $X $ is homeomorphic to the (standard) cone over the link,}
  $Cone[Link[X]]\sset (\C^N,o).$
   In ``most cases" this (non-embedded) homeomorphism cannot be chosen differentiable in any sense.

  In some cases this  topological equivalence can be strengthened to (non-embedded) bi-Lipschitz equivalence, as follows.
   The (standard) metric on $(\C^N,o)$ induces the inner metrics on $X $ and $Cone[Link[X]],$ by the length of the shortest path between two points.
    The germ $X $ is called inner-metrically-conical (IMC) if the homeomorphism $X \isom{} Cone[Link[X]]$ can be chosen bi-Lipschitz \wrt
     these inner metrics.
 The metrically conical structure is the simplest possible metric structure. E.g. if $X$ metrically conical, then it is bi-Lipschitz equivalent to   its triangulation.

Any complex  analytic curve germ is IMC. In higher dimensions the situation is more complicated.
 Initially many complex-analytic surface germs were thought to be IMC. The first
  weighted-homogeneous non-IMC surface germs were found in \cite{Birbrair-Fernandes.08}.
   If $X $ is a weighted-homogeneous IMC-germ of quotient type, $\quots{\C^2}{\mu},$ then its two lowest weights are necessarily equal, \cite{Birbrair-Fernandes.Neumann.08}.
   The converse statement was verified in \cite{Birbrair-Fernandes.Neumann.09} for  Brieskorn-Pham singularities in $(\C^3,o)$.
 Finally, \cite{Birbrair.Neumann.Pichon.14} established the thick-thin decomposition for normal surface germs.
 In particular, this gave an algorithm to verify the (non)IMC-property via the resolution graph of the singularity.
  See also \cite{Birbrair.Fernandes.Grandjean.17}.
  Using this \cite{Okuma.17} proved: a Brieskorn-Pham surface germ that is an isolated complete intersection is IMC iff its two lowest weights are equal. 
 In higher dimensions, among the $A_k$-types (with equation $\sum^n_{i=1} x^2_i+x^{k+1}_{n+1}=0$), the only IMC-case is $A_1,$ \cite{Birbrair.Fernandes.Neumann.Grandjean.O'Shea14}.

In all these cases the first reason be non-IMC is the existence of so-called ``fast loops'', and their generalizations, ``choking horns", \cite{Birbrair.Fernandes.Grandjean.17} These are topologically  non-trivial loops on $Link[X]$ that vanish  faster than linearly, as one  approaches the base point $o\in X.$

  Fast loops (and fast cycles in higher dimensions) are important far beyond the IMC-context.
   They carry   essential information on the (Lipschitz part of) local geometry and topology of $X .$
    E.g., each fast loop lies inside a thin zone, a rigid sector of the germ, where the metric is essentially distorted (as compared to the metric on  a cone).
To detect/forbid these fast loops is in general a complicated task.
 The only well-studied case is the normal surface singularities,  \cite{Birbrair.Neumann.Pichon.14}, with the fast-loops criteria via the resolution graph of $X$.

 In \cite{Kerner-Mendes.Weighted.Homogen} we gave sufficient criteria for fast cycles on complete intersection germs of any (co)dimension  that are perturbations of weighted-homogeneous germs. In ``most cases" if the $n$ lowest weights do not all coincide, then $X$ has a fast cycle.

All these results suggested:   IMC-germs are very rare, to the extent that, e.g. complex-analytic IMC's in $(\C^3,o)$ could be possibly classified, up to the ambient topological equivalence.

 Now we give simple necessary/sufficient conditions  for reduced germs to posses fast loops. These germs are   of any dimension and codimension, not necessarily normal or complete intersection, possibly with non-isolated singularities.
\bei
\item For surface germs  the conditions are in terms of the discriminant of the covering. They are necessary and sufficient.
\\For \!normal \!surface \!germs \!this \!allows \!simple \!verification \!of \!the \!IMC \!property.
 \!In \!particular \!one gets: while IMC's are rare,
  their classification is in some sense impossible. (Even in the case $X \sset (\C^3,o)$.)

\item
The higher dimensional case admits (in many cases) a reduction to the surface case. Hence an easy way to detect a fast cycle. (And thus an obstruction to the IMC-property.)

\eei

Some of our results (for normal surface germs) could be possibly obtained by the methods of \cite{Birbrair.Neumann.Pichon.14}. But our proofs do not use resolution of singularities and are independent of those of  \cite{Birbrair.Neumann.Pichon.14}.

\

A remark: our results allow to detect the homotopy types of fast loops, and their vanishing speeds. This will be done in the subsequent work.

\subsection{The fast cycles and (non)IMC-criteria}
\subsubsection{}\label{Sec.Intro.Fast.Cycles.Surfaces}
Let $X $ be a reduced, complex-analytic surface germ,  of  multiplicity $p,$ not necessarily normal, possibly with non-isolated singularity. Present it as a covering space,  $\pi:X \stackrel{p:1}{\to}(\C^2,o).$ (This covering is assumed to be ``convenient".)
 This map is ramified over the discriminant, $\De\sset (\C^2,o).$
 We detect fast loops by the explicit numerical condition.
\\{\bf Theorem \ref{Thm.Fast.Loop.criterion.discriminant.surfaces}.} {\em $X $ has no fast loops iff the inequality $p> r_k-1+\sum q_i\cdot mult(\De_{k,i})$ holds for each non-smooth tangential component of $\De.$}

The ingredients $\{r_k\},$  $\{q_{i}\},$  $\{mult(\De_{k,i})\}$ form the ramification data of the covering $X \stackrel{p:1}{\to}(\C^2,o),$ and are defined in \S\ref{Sec.Covering.Criteria.Ramification.Data}. In many cases this data is easily computable from the equations of $X.$ Thus the theorem is handy in applications.
  In fact, as the proof shows,
 $r_k +\sum q_i\cdot mult(\De_{k,i})-p$ is the lower bound on the number of fast loops ``hanging over $\De_k$".
\\{\bf Corollary \ref{Thm.IMC.criterion.Corollary}.} {\em A normal surface germ is IMC iff
 $p> r_k-1+\sum q_i\cdot mult(\De_{k,i})$ for each non-smooth tangential component of $\De.$}

The absence of fast loops, i.e. this numerical inequality, imposes strong restrictions on the projectivized tangent cone.
Let $X\sset (\C^{2+c},o) $ be an isolated strictly complete intersection surface germ.
 Its  projectivized tangent cone, $\P T_{X }\sset \P^{1+c},$ is a curve  (singular,  possibly non-reduced).
  Take its degree, $deg[\P T_X],$ and the degree of is reduced locus, $deg[\P T^{red}_X].$
 Suppose the germs $(\P T_X,pt_i)$  are reduced and singular for some points $\{pt_i\}.$ Take their Milnor numbers,
  $\mu_i:=\mu(\P T_{X },pt_i).$
\\{\bf Theorem \ref{Thm.Fast.Loops.via.Tangent.Cone.surfaces}.}
 {\em 1. Suppose these points $\{pt_i\}$ lie in one hyperplane, $\P^c\sset \P^{1+c},$ that is ``non-tangent" to $\P T_X.$
  If $\sum \mu_i \ge 1+deg[ \P T_{X }]-deg[ \P T^{red}_{X }],$ then $X $ has a fast loop.
\\2. In particular, if  $\P T_{X }$ is reduced, and $X $ is   IMC, then $\P T_{X }$ must be smooth. (Hence $X $ is an ordinary multiple point.) Additionally, when $\P T_{X }=\P T_X^{red}$ any reduced singular point of $\P T_X\cap \P^c$ produces a fast loop.}

  This is a very simple (and useful) IMC-obstruction. For example:
   \bei
   \item The only IMC's among the super-isolated surface germs are the ordinary multiple points.
    (Example \ref{Ex.IMC.conditions.imposed.by.PT_X})
\item The only IMC's among the surface singularities in $(\C^3,o)$ of right modality$\le2$  are: $A_1,D_4,X_9,P_8.$
 (Corollary \ref{Thm.IMC.among.right.modality.le.2})
\eei

  \subsubsection{}
 In   higher-dimensional case our results are weaker, giving just a sufficient condition for fast loops.  Take a convenient covering $\pi:X \to (\C^n,o),$ with the discriminant $\De\sset (\C^n,o),$ and restrict it to a smooth surface germ, $\pi_\cS: X \cap \pi^{-1}\cS\to \cS\cong (\C^2,o).$ Accordingly one gets the discriminantal curve $\De_\cS:=\De\cap \cS,$ and the ramification data,
   $ r_k ,$  $\{q_{i}\},$  $\{mult(\De_{k,i})\},$ as above. Suppose $\cS$ is ``weakly non-special", i.e. the same numerical data is obtained for any other surface that is {\em tangent} to $\cS.$
\\{\bf Proposition \ref{Thm.Fast.Loops.via.discriminant.higher.dim}.} {\em If $p\le r_k-1+\sum q_i\cdot mult(\De_{k,i})$   for some non-smooth tangential component of $\De_\cS$ then $X $ has a fast loop.}

This gives a quick check of fast loops/non-IMC property in higher dimensions. E.g. (corollary \ref{Thm.Fast.Loops.via.PT.discriminant.higher.dim}): if $\P T_{\De}$ has a point of high multiplicity, then $X$ has a fast loop.

\subsection{Acknowledgements} Our thanks are to Lev Birbrair, Alexandre Fernandes, Andrey Gabrielov, Edson Sampaio, for important advices.

\subsection{The contents/structure of the paper}\label{Sec.Introduction.Contents.of.the.paper}

\bee[\!\S1\!] \setcounter{enumi}{1}
\item
 contains basic definitions/facts on singular germs. They are used in \S3 to detect fast loops.
 \bee[\!\S2.1\!]\setcounter{enumii}{1}
\item
recalls the  tangent cone, $T_{X },$ and its projectivization, $ \P T_{X },$ for complete intersections.
 In the later case $\P T_{X }$ is a complete intersection only if $X $ is a ``strictly complete intersection".
\item
defines ``convenient" coverings $\pi:X \to (\C^n,o).$  The critical/discriminantal loci are defined set-theoretically.
\item describes the  blowup of the covering $\pi$, and the finiteness/flattness of the strict transform $\tpi.$

\item recalls the standard scheme structure on the critical locus when  $X $  is a complete intersection. Then the critical locus is also a complete intersection. Its image (with the Fitting scheme-structure) is the discriminant, a hypersurface.

   \item The covering $\pi:X \to(\C^n,o)$ is a Lipschitz map. It is far from being locally bi-Lipschitz.
    Yet, we have: the preimage of a thin set is thin. In particular, any pair of tangent arcs in $(\C^n,o)$ can be lifted to a pair of inner-tangent arcs in $X .$
\item recalls fast cycles, the basic obstruction to the IMC-property.
 For normal surface germs the fast loops are the only obstructions. For non-normal surface germs there are various other obstructions.
\eee
\item detects fast loops on $X$.
\bee[\S3.1]
\item defines the covering data of $\pi:  X\to (\C^2,o)\supset \De,$ i.e. the tangential components $\{\De_k\}$ of $\De,$ the
 total ramification orders $\{q_i\},$ and the special invariant $r_k.$
\item contains Theorem \ref{Thm.Fast.Loop.criterion.discriminant.surfaces}: the condition
 $p\le r_k-1+\sum q_i\cdot mult(\De_{k,i})$ means the existence of a fast loop.

\hspace{0.3cm}\parbox{15.3cm}{In the proof one studies the covering $X\to (\C^2,o).$ Suppose the discriminant $\De$ has a non-smooth tangential
component, $\De_k\sset \C^2_{x_1,x_2}.$
 Let its tangent line be $Span_\C(\hx_1).$ The intersection $\De_k\cap \{Im(x_1)=0\}$ consists of several arcs. Connecting them inside the planes  $\{Re(x_1)=t\}$
 produces a (subanalytic) H\"{o}lder polyhedron.
 (Which is  a thin set.)
  Its lifting, $H\ddot{o}ld^X\sset X,$ is still a thin set, by \S\ref{Sec.Preliminaries.Thin.Sets}.
    It remains to verify: one of the components of $H\ddot{o}ld^X\smin o$ is not ``linkwise contractible". Which means: for $\ep>0$ small enough, one of the connected components of
     $H\ddot{o}ld^X\cap V(x_1-\ep)$ is non-contractible. And this is verified via the homology $H_1(H\ddot{o}ld^X\cap V(x_1-\ep)),$
     by computing the Euler characteristic.}

\

 Then go explicit examples. We classify all the IMC-surface germs of multiplicity two. Then we classify   IMC's
 among the surfaces of type $V(z^p+a_1(x,y)z+a_o(x,y))\sset (\C^3,o).$

\item gives the condition  on the tangent cone imposed by fast-loop, Theorem \ref{Thm.Fast.Loops.via.Tangent.Cone.surfaces}.
 The proof uses Theorem \ref{Thm.Fast.Loop.criterion.discriminant.surfaces} and the facts established in
\S\ref{Sec.Preliminaries.Blowup.of.pi}, \S\ref{Sec.Preliminaries.Crit.And.Discr}.

Now fast loops are detectable immediately and we deduce the applications mentioned in \S\ref{Sec.Intro.Fast.Cycles.Surfaces}.

 \item is  the higher dimensional case. We give a sufficient condition for fast loops via the ramification data
  (proposition \ref{Thm.Fast.Loops.via.discriminant.higher.dim}). In some cases even singularities of $\P T_\De$ that have high multiplicity
  ensure  vanishing cycles (corollary \ref{Thm.Fast.Loops.via.PT.discriminant.higher.dim}).
\eee

\eee

\subsection{Notations and conventions}\label{Sec.Preliminaries.Dictionary}
\bee[\bf i.]

\item Through the paper by germs we mean their small representatives.

Everywhere in the paper $X\sset(\C^N,o)$ is a reduced pure-dimensional complex-analytic germ, of dimension $n\ge2,$ and multiplicity $p,$ possibly with non-isolated singularity.

All the other germs are  subanalytic and  closed. 

\item All our arcs are (real) subanalytic. Usually we take their length-parametrization, $\ga(t)=(s_1\cdot t+o(t),\dots,s_N\cdot  t+o(t)),$ with $(s_1,\dots,s_N)\neq (0,\dots,0).$

By a foliation we always mean a singular foliation of the germ $X$
  by real arcs, all passing through the origin.

\item The (inner) tangency order of two arcs $\ga_1, \ga_2 \sset X$ is $tord_X (\ga_1,\ga_2):=ord_t d_X(\ga_1(t),\ga_2(t)).$
 Here $d_X:   X\times X\rightarrow \R$ is a subanalytic distance equivalent to the distance given by the infimum of lengths of paths connecting points of $X$ \cite{Kurdyka.Orro.97}. Note that $tord_X (\ga_1,\ga_2)$ can be non-integer.

 If $X=(\R^N,o)$ then $d_X=\|\dots\|$ (the usual Euclidean distance) and we write just $tord(\ga_1,\ga_2).$

 By a LNE germ $X$ means a representative where the inner distance $d_X$ and the outer distance $\|\dots\|$ on $X$ are equivalent.

The tangency order of an arc and a subgerm,  $\ga , Y\!\sset\! X,$ is $tord_X(\ga,Y\!)\!:= \!\supl_{\!\ga_Y\sset Y} tord_X(\ga,\ga_Y)\!\le\!\infty.$
The ``order of difference" of
two subgerms,   $Y_1 , Y_2\!\sset\! X,$ is
\beq\label{Eq.Order.of.difference}
 dif\!f.ord_X(Y_1,Y_2):=\infl_{\ga_1\sset Y_1} tord_X(\ga_1,Y_2) =\infl_{\ga_2\sset Y_2} tord_X(\ga_2,Y_1) \le \infty.
 \eeq

\item
 Take  an $\R$-analytic germ   $Z\sset (\R^N,o) .$ The tangent cone $T_Z$ can be defined set-theoretically, e.g. as the Whitney cone $C_3$,
  \cite[Chapter 7]{Whitney},  \cite[Chapter 2]{Chirka}. The definition of $T_Z$ is the same for subanalytic germs.
  The tangent cone is subanalytic as well, with $dim_\R T_Z \le dim_\R Z.$ In particular,
    if  $dim_\R T_Z=1,$ then $T_Z$ is a finite union of lines and half-lines.

\item Let $Y\sset \R^N$ be a sub-analytic set.
Every homology class in $H_*(Y,\Z)$ admits a subanalytic representative. Indeed, the singular homology is isomorphic to the simplicial homology.
 And sub-analytic manifolds admit subanalytic triangulations, \cite[\S4.3]{Coste}.

\item The singular locus of a   complete intersection $X=V(f_1,\dots,f_c)\sset (\C^{n+c},o)$ is defined (inside $X$) by the  ideal of maximal minors $I_c[f'_1,\dots,f'_c].$ Here $\{f'_j\}$ are the column gradients.

\eee

\section{Preparations}

\subsection{The (projectivized) tangent cone for complete intersections}\label{Sec.Preliminaries.Tangent.Cone}
\subsubsection{Hypersurfaces}
 Let $X\! \!=\!V(f)\!\sset \!(\C^N,o),$ then $T_{X }\!=\!V(l.t.(f))\!\sset\! \C^N,$  the lowest-order terms of $f.$
    Observe that the polynomial $l.t.(f)$ can be not square-free, in which case $T_{X }$ is non-reduced.

The projectivized tangent cone is  $\P T_{X }=V(l.t.(f))\sset \P^{N-1}.$ This is a projective hypersurface, possibly reducible, non-reduced.
  Its (total) degree   equals the multiplicity of $X ,$ i.e. the Taylor order $ord(f).$

\subsubsection{Complete intersections}\label{Sec.Preliminaries.Tangent.Cone.C.I.}
 Let $X =V(\uf )\sset (\C^{n+c},o),$ where $\uf:=f_1,\dots, f_c\in \C\{\ux\}$ is a regular sequence.
  Then $T_{X }=V(l.t.(\uf)),$ defined by the lowest order terms of all the elements of the ideal $(\uf).$
   In particular, $l.t.(\uf)\supseteq (l.t.( f_1),\dots, l.t.(f_c)).$ This inclusion can be proper.
 \bex
 For $f_1(x,y)\!=\!x^2y\!+\!y^4,$  $f_2(x,y)\!=\!xy^2\!+\!x^4,$ one has $l.t.(\uf)\!=\!(x^2y,xy^2,x^5\!-\!y^5)\!\supsetneq \!(l.t.( f_1),l.t.( f_2)).$
\eex
Therefore $T_{X }\sset \C^{n+c}$ is not necessarily a complete intersection.
 The germ $X $ is called a {\em strictly complete intersection} if both $X $ and $T_{X }$ are complete intersections, \cite{Benett-1977}.
 In this case:
 \bei
 \item There exists a choice of defining equations, $X=V(f_1,\dots,f_c)$ with $T_X=V(l.t.(f_1),\dots,l.t.(f_c))$ and
  $deg(\P T_{X })=multX=\prod ord(f_i) .$
 \item The projectivization $\P T_{X }\sset \P^{n+c-1}$ is a complete intersection (possibly reducible, non-reduced).
 \item The strict transform under blowup, $\widetilde{X}\sset Bl_o(\C^{n+c},o),$ is a locally complete intersection at all its points.
\eei

\subsection{The convenient covering and its critical locus}\label{Sec.Preliminaries.Convenient.Coverings}
Let $X  $ be as in \S\ref{Sec.Preliminaries.Dictionary}.i.
 Take the coordinate plane $\C^n= V(x_{n+1},\dots,x_N)\sset \C^N$ and the projection $\C^N\ni\ x\to (x_1,\dots,x_n)\in\C^n$.
 We call the   restriction of  this projection to the germ, $X \stackrel{\pi}{\to}(\C^n,o)$,  {\em convenient} if:

 \bee[\bf i.]
 \item $\pi$ is a (finite) ramified  covering;
 \item The kernel of the projection  at $o$ is non-tangent to $X ,$ i.e.
 $T_X\cap V(x_1,\dots,x_n)=o\in \C^N.$
 \eee
 These assumptions hold for generic projections (observe that the condition ii. is enough to get convenience). A projection can be quite non-generic, and yet convenient. E.g. the covering $(\C^3,o)\supset V(z^p-f(x,y))\to (\C^2,o)$ is totally ramified over the
\\\parbox{12cm}
{curve germ $V(f)\sset (\C^2,o).$ This covering is convenient iff $ord\ f\ge p.$ The critical locus of $\pi$ consists of those points where $\pi$  is locally not an isomorphism \!onto \!its \!image.
 \!Thus \!$Crit(\pi)\!\sset\! X $ \!is \!an \!analytic \!subgerm, \!and
 $dim_\C Crit(\pi)\!\le\! (n-1).$ \!(By Sard's theorem.)
}\hspace{1cm} $\bM (\C^N,o) \supset X  \supset Crit(\pi)
  \\\quad \quad \quad  \searrow\quad \downarrow\pi \quad\quad\downarrow\pi_|
  \\\quad\quad\quad\quad (\C^n,o) \supset\De\quad\eM$
\\\!We remark: \!$Crit(\pi)\!\supseteq\! Sing(X).$

The discriminant of $\pi$  is the image of the critical locus, $\De:=\pi(Crit(\pi))$. As $\pi$ is a finite
 morphism, $\De\sset (\C^n,o)$ is an analytic subgerm, $dim_\C\De\le n-1.$ (By Grauert's direct image theorem.)
 The inequality can occur here. E.g. let $X $ be the union of two planes (of dimension $n$) intersecting at one point.
  Then $\De$ is just one point, $o\in\C^n.$

 When $X $ is a complete intersection,   these set-theoretic definitions are replaced by more precise algebraic ones,
 \S\ref{Sec.Preliminaries.Crit.And.Discr}, and $\De$ is well-behaved, e.g. $dim_\C\De= n-1.$

Recall the basic properties of convenient projections.
\bei
\item
 A small deformation of the projection $\pi$ preserves  the convenience of the covering.
 \item
 If $\pi$ is convenient, then  the restriction $\pi_Y=\pi|_{\pi^{-1}Y}$ is a convenient covering for any smooth subvariety $Y\sset (\C^n,o)$, assuming $Y\not\sseteq \De$ and $\pi^{-1}Y$ is reduced.
 
\item A small deformation of the base point does not preserve the convenience. E.g., let $o\neq x_o\in Crit(\pi)$ be a smooth point of $X.$
 Then the projection $(X,x_o)\to (\C^n,\pi(x_o))$ is non-convenient.
\eei

\subsection{The blowup of $\pi$}\label{Sec.Preliminaries.Blowup.of.pi}
Let $X  $ be as in \S\ref{Sec.Preliminaries.Dictionary}.i.
 Suppose the projection $\pi\colon X \to(\C^n,o)$ is  convenient.

\bel \label{Thm.Diagram.Blowup.Lift.Covering}
\bee
\item
The  blowup of the origin produces the diagram:
\[
\bM  (\C^N,o)&\supset  X &\leftarrow& \tX:=Bl_oX &\supset&\P T_{X }=\tX\cap \P^{N-1}
\\   &\pi \downarrow&&\tpi\downarrow&&\downarrow \tpi_|&& \\
  &(\C^n,o)&\leftarrow &  Bl_o(\C^n,o)&\supset& \P^{n-1}\eM
 \]
Here the morphisms $\tpi,\tpi_|$ are finite.
\item
 If $X \!\sset \!(\C^N,o)$ is a  strictly complete intersection    (\S\ref{Sec.Preliminaries.Tangent.Cone.C.I.}), then  the morphisms $\tpi,\tpi_|$ are flat.
\eee
\eel
\bpr
\bee
\item
The  projection $Bl_o\C^N\dashrightarrow Bl_o \C^n$ is $(x_1,\dots,x_N,[\si_1: \dots :\si_N])\to (x_1,\dots,x_n,[\si_1:\dots:\si_n])$, it is defined
 only  for $(\si_1,\dots,\si_n)\neq (0,\dots,0)$. But $\pi$ is a convenient projection, therefore $\tX\cap V(\si_1,\dots,\si_n)=\empty$.
  Therefore the projection $\tpi: \tX\to Bl_o(\C^n,o)$ is well defined.

All the fibres of $\tX\smin \P T_X\to \C^n\smin \{o\}$ are finite.
 And the restriction of the projection to the exceptional divisor, $\tpi_|:\tX\cap \P^{N-1}\stackrel{}{\to}\P^{n-1}$, has finite fibres.
 Indeed, if a fibre over a point $[\si_1:\dots:\si_n]\in \P^{n-1}$
  is infinite, then  $ \tilde\pi^{-1}(0,\dots,0,[\si_1:\dots:\si_n])\sset \P T_{X }\sset \tX$ is  an  algebraic subset of positive dimension. But then $dim(T_{X }\cap V(x_1,\dots,x_n))>0$,
 contradicting the convenience of $\pi$.
\\Thus $\tpi $ is a projective morphism with finite fibres. Hence $\tpi$ finite. And therefore $\tpi_|$ is finite as well.
 \item
 As $X $ is a strictly complete intersection, both $\tX$ and $\P T_{X }$ are   locally complete intersections, while $\tpi$ and $\tpi_|$ are finite morphisms with smooth images.

  Thus to verify the flatness it is enough to verify: the degree of fibres is constant.
   For each point $\tx\in Bl_o(\C^n,o)\smin \P^{n-1}$ the fibre $\tpi^{-1}(\tx)$ is  of degree $p$ (=the degree of the covering $\pi$).
    For each point  $\tx\in  \P^{n-1}$ we have:
  \beq
  deg \tpi^{-1}(\tx)=  deg\ \tpi_|^{-1}(\tx)=deg\ (\P T_{X } \cap \P^{N-n})\stackrel{*}{=}deg (\P T_{X } )\stackrel{**}{=}multX=p .
  \eeq
 Here (after a $GL(\C^n)$ transformation) $\P^{N-n}=V(\si_1,\dots,\si_{n-1})\sset \P^{N-1}.$
  The equalities  $*$,  $**$ hold  because  $X,\P T_{X }$  have no embedded components of  dimension$<n$ (being  complete intersections).
\epr
 \eee

\subsection{The scheme structure on the critical locus and the discriminant}\label{Sec.Preliminaries.Crit.And.Discr}
Let $X$ be as in \S\ref{Sec.Preliminaries.Dictionary}.i. Let
$X \stackrel{\pi}{\to}(\C^n,o)$ be a (not necessarily convenient)  ramified $p:1$ covering.
\bel\label{Thm.Discriminant.Basic.Properties}
If the germ $X \sset (\C^{n+c},o)$ is a complete intersection then the critical
 locus and the discriminant (with the scheme structure defined in the proof) have the following properties.
\bee
\item $Crit(\pi)\sset (\C^{n+c},o)$ is a complete intersection of dimension $n-1$,   $\De\sset (\C^n,o)$ is a hypersurface.

\item Let $X \!=\!V(\uf)\!\sset \!(\C^{n+c},o),$ for
    $\uf=(f_1,\dots,f_c),$ the column.
 Denote $p_i: =\!ord(f_i).$ Then

 \hspace{2cm} $mult(Crit(\pi))\!\ge\!(\prod p_i)\cdot\! ord (det[  \uf']) \ge\!(\prod p_i)\cdot  \sum(p_j-1).$
\\
Here $[  f']\in Mat_{c\times c}$ is the Jacobian matrix of derivatives with respect to the variables $x_{n+1}\dots x_{n+c}.$
\item The diagram of \S\ref{Sec.Preliminaries.Convenient.Coverings} is functorial on pullbacks.
 Namely, for each morphism $\phi:(\C^m,o)\to (\C^n,o) $ with $\phi(\C^m,o)\not\sseteq \De,$
 and the corresponding projection $\phi^*\pi: X \times_{(\C^n,o)}(\C^m\!,o)\to(\C^m\!,o)$
one has:
    \bei
    \item $Crit(\phi^*\pi)=\phi^*Crit(\pi) ,$ i.e. $\cO_{Crit(\phi^*\pi)}=\cO_{Crit(\pi)}\otimes_{(\C^n,o)}\cO_{(\C^m,o)}$,
     \item  $\De_{\phi^*\pi}=\phi^*\De ,$ i.e. $\cO_{\De_{\phi^*\pi}}=\cO_\De\otimes_{(\C^n,o)}\cO_{(\C^m,o)} .$
\eei
 \eee
\eel
Here $X \times_{(\C^n,o)}(\C^m\!,o):=Spec(\cO_{X } \otimes_{  (\C^n,o) } \cO_{(\C^m,o)})$ is the pullback of $X $ via $\phi.$
\bpr First we define $Crit(\pi)\sset X$ and $\De\sset (\C^n,o)$ as analytic subgerms. Then we deduce the statements.
\bee[\bf a.]
\item
    Recall the standard definitions  \cite{Looijenga}, \cite{AGLV.II}. When $X$ is smooth, the critical locus of the projection $X \to (\C^n,o)$ is the degeneracy locus of the induced morphism
     $T_{X }\to T_{(\C^n,o)}.$ In the general case the  critical locus is defined as the degeneracy locus of the $n$-form $dx_1\wedge \cdots\wedge dx_n|_{X }\in \Om^n_{X }$.
  For the defining equations, $ X  =V(f_1,\dots,f_c)\sset (\C^{n+c},o),$ one has  $\Om^n_{X }=\cO_{X }\otimes\quot{\Om^n_{(\C^{n+c},o)}}{Span\{df_i\}\wedge \Om^{n-1}_{(\C^{n+c},o)}}.$
    Thus the critical locus is the set of points of $X$ where the form $dx_1\wedge \cdots\wedge dx_n|_X$ belongs   to
  the submodule  $Span\{df_i\}\wedge \Om^{n-1}_{(\C^{n+c},o)}|_{X }$. Equivalently,
     $dx_1\wedge \cdots\wedge dx_n\wedge df_1\wedge\cdots\wedge df_c|_{X }=0.$ Writing this in coordinates, the critical locus is defined (inside  $X$)
      by the equation $det[\di_{n+1}\uf\dots \di_{n+c}\uf]=0.$

  \noindent    We get the analytic subgerm $Crit(\pi)\sset (\C^n,o)$ and its local ring
  \beq\label{Eq.local.ringof.Crit}
  \cO_{Crit(\pi)}=\quot{\cO_{(\C^{n+c},o)}}{(\uf,det[ \di_{n+1}\uf,\dots,\di_{n+c}\uf])}.
  \eeq
\item
We claim: the subgerm $Crit(\pi)\sset (\C^{n+c},o)$ is a  complete intersection, i.e. the sequence $f_1,\dots,f_c,$ $det[\dots]$ is regular in
 $\cO_{(\C^{n+c},o)}.$
        Indeed, $dim(Crit(\pi))\ge n-1,$ as the critical locus is defined by $(c+1)$ equations. But
        the
      projection $X \to (\C^n,o)$ is generically unramified on every irreducible component of $X.$
      Therefore $dim(Crit(\pi))\le n-1.$ Altogether, the sequence $f_1,\dots,f_c,det[\dots]$ of length $(c+1)$ defines a subgerm $Crit\sset(\C^{n+c},o)$ of codimension $(c+1).$ Hence this sequence is regular, \cite[\S17]{Eisenbud-book}.

Equation \eqref{Eq.local.ringof.Crit} gives then Part 2 of the statement.

 The discriminant is the image of $Crit(\pi).$  As $\pi$ is a finite morphism, $\De$ is an analytic subset,  $dim(\De)=n-1,$ i.e. $\De\sset (\C^n,o)$ is
      (set-theoretically) a hypersurface.

 The natural scheme structure on $\De$ is the image structure.  The pushforward $\pi_*\cO_{Crit(\pi)}$ is a module over $\cO_{(\C^n,o)}.$ Its support is $\De\sset (\C^n,o)$, of dimension $n-1.$  Therefore we define $\De\sset (\C^n,o)$ as an analytic subgerm by the zeroth
      Fitting ideal, $Fitt_0[\pi_*\cO_{Crit(\pi)}]\sset \cO_{(\C^n,o)},$ \cite[\S20]{Eisenbud-book}.
 In detail, take a presentation matrix, $\pi_*\cO_{Crit(\pi)}=coker[A].$ Below we show that $A$ is a square matrix. Then $Fitt_0[\pi_*\cO_{Crit(\pi)}]=(det[A]).$

 \item
As $\De\sset (\C^n,o)$ is a (set-theoretic) hypersurface germ, we can assume: the intersection $V(x_1,\dots,x_{n-1})\cap \De$ is a point.
 Thus $x_1,\dots,x_{n-1}$ goes to a regular sequence in the ring $\cO_\De.$

 We claim: the sequence $x_1,\dots,x_{n-1}$ is regular also on the ring $\cO_{Crit(\pi)},$
  and hence on the module  $\pi_*\cO_{Crit(\pi)}$ as well.
  First, $x_1$ is not a zero-divisor on $\cO_{Crit(\pi)}.$
 Indeed,  $Crit(\pi) $ has no embedded components, being a complete intersection. In addition, $V(x_1)$ intersects  $Crit(\pi)$ properly.

  Now restrict onto $V(x_1).$ One has:
\beq
\cO_{Crit(\pi)}\otimes \quot{\cO_{(\C^{n+c},o)}}{(x_1)}\!=\!\quot{\cO_{(\C^{n+c},o)}}{(x_1,\uf ,det[\dots])}\!\cong\!
\quot{\cO_{(\C^{n+c-1},o)}}{(\uf|_{V(x_1)} ,det[\dots|_{V(x_1)}])}\!=\!\cO_{Crit(\pi|_{V(x_1)})}.
\eeq

The restricted map $X|_{V(x_1)}\to (\C^{n-1},o)$ is still a ramified covering. Therefore the restrictions
 $(f_1|_{V(x_1)},\dots,f_c|_{V(x_1)},det[\dots]|_{V(x_1)})$ still form a regular sequence and $\cO_{Crit(\pi|_{V(x_1)})}$ is a complete intersection.
 Repeat the argument for $x_2$, and so on.

\item

    The sequence $x_1\dots x_{n-1}$ is regular on the module $\pi_*\cO_{Crit(\pi)}.$
 Therefore the $\cO_{(\C^n,o)}$-module $\pi_*\cO_{Crit(\pi)}$  is of $depth=n-1=dim(\De),$ \cite[\S18]{Eisenbud-book}.
   Hence $\pi_*\cO_{Crit(\pi)}$ is Cohen-Macaulay.
  By Auslander-Buchsbaum formula, the $\cO_{(\C^n,o)}$-resolution of $\pi_*\cO_{Crit(\pi)}$ (over $\cO_{(\C^n,o)}$) is of length one.
 Thus the presentation matrix,  $\pi_*\cO_{Crit(\pi)}=coker[A],$ is square. Hence $Fitt_0(\pi_*\cO_{Crit(\pi)})=(det[A])$ is a  principal ideal.
 Thus the subscheme $\De\sset (\C^n,o)$ is a  hypersurface (possibly with multiple components).
\epr
\eee
\bex\bee[\bf i.]

\item
Let $X :=V(z^p-f(x,y))\sset (\C^3,o),$ where $ord(f), p\ge2.$  Take the covering $X\to (\C^2_{xy},o),$ it is convenient iff $ord(f)\ge p.$
 Then the defining ideal of $Crit(\pi)\sset (\C^3,o)$ is $(z^{p-1},f(x,y))\sset \cO_{(\C^3,o)}.$
 Set-theoretically one gets: $\De=V(f(x,y))\sset (\C^2,o).$ To obtain the image scheme
 structure on $\De$ we consider $\cO_{Crit(\pi)}=\quots{\C\{x,y,z\}}{(z^{p-1},f(x,y))}$ as an $\C\{x,y\}$-module.
  Its natural generators are $\{1,z,\dots,z^{p-2}\},$ and the module splits as a direct sum.
  We get the presentation
\beq
\oplusl^{ p-1} \C\{x,y\} \stackrel{\oplus[f ]}{\longrightarrow}\oplusl^{ p-1} \C\{x,y\} \to \cO_{Crit(\pi)}\to0.
\eeq
   Hence the defining ideal of the discriminant is $I_\De=det[\oplus[f]]=(f)^{p-1}\sset \C\{x,y\}.$

\item (Functoriality on restrictions) Take a linear subspace $(\C^m,o)\sset(\C^n,o)$ with $(\C^m,o)\not\sseteq \De.$ Take the restricted projection
$\pi_|:\ \pi^{-1}(\C^m,o)\to (\C^m,o).$ One gets: $Crit(\pi_|)=Crit(\pi)\cap \pi^{-1}(\C^m,o)$
 and $\De(\pi_|)=\De\cap (\C^m,o).$

\eee\eex

\subsubsection{The discriminant for hypersurface germs}
Take a hypersurface germ $X=V(f)\sset (\C^{n+1},o)$   of multiplicity $p$. Suppose the projection $X\to (\C^{n}_{x_1,\dots,x_n},o)$ is convenient.
 Thus $T_X\cap V(x_1,\dots,x_n)=(o),$
 therefore $f$ contains the monomial $x^p_{n+1}.$ By Weierstra\ss\ preparation theorem one can present:
\beq
f(x)=u(x)\big(x^p_{n+1}+x^{p-1}_{n+1}\cdot a_{p-1}(x_1,\dots,x_n)+\cdots+a_0(x_1,\dots,x_n)\big).
\eeq
 Here $u$ is invertible and $a_0,\dots,a_{p-1}$ are analytic power series, with  $ord(a_j)\ge p-j$.
  Thus we omit $u.$


This hypersurface is the pullback of the universal hypersurface, $V(z^p+z^{p-1}a_{p-1}+\cdots+a_0)\sset \C^{1+p}_{z,a_{p-1},\dots,a_0},$
 under the map $(x_1,\dots,x_n)\to (a_{p-1}(x),\dots,a_0(x)).$
Therefore (by part 3 of lemma \ref{Thm.Discriminant.Basic.Properties}) the discriminant of the projection, $\De\sset \C^n,$   is obtained as the pullback of the classical discriminant $\De_{class}$ of the polynomial
 $z^p+z^{p-1}a_{p-1}+\cdots + a_0$. Here $\De_{class}\in \Z[\{a_\bullet\}]$ is a (complicated) polynomial in the coefficients $\{a_\bullet\}$. We use the following properties \cite{GKZ}:
 \bei
 \item Assign  the weights to the coefficients, $w(a_j)=p-j$, then the polynomial  $\De_{class}$ is  weighted-homogeneous, of total weight $p(p-1)$.
 \item $\De_{class}(\{a_\bullet\})=(\frac{a_0}{p-1})^{p-1}-(\frac{a_1}{p})^p+ g(\{a_\bullet\})$,
  where $g(\{a_\bullet\})\in (a_0,a_1)\cdot (a_2\dots a_{p-1})$.
\item For $p\ge3$ the hypersurface $V(\De_{class})  \sset \C^p$ is singular in codimension 1.
 \eei

\subsubsection{The polar multiplicity of a curve germ}\label{Sec.Preliminaries.Crit.and.Discr.Polar.Multipl}
   Take a reduced complete intersection curve germ, $(C,o)=V(f_1,\dots,f_c)\sset (\C^{1+c},o).$ Suppose the projection
    $(\C^{1+c},o)\supset (C,o)\to (\C^1_{x_{1+c}},o)$ is convenient, \S\ref{Sec.Preliminaries.Convenient.Coverings}.
      The critical locus of the projection is $V(f_1,\dots,f_c)\cap V(det[\di_i f_j]_{i,j\le c}),$ see equation \eqref{Eq.local.ringof.Crit}.
       This is a one-point scheme, a complete intersection. Its degree is  the polar multiplicity,
        $m_1(C,o)=deg[V(f_1,\dots,f_c)\cap V(det[\di_i f_j]_{i,j\le c})],$ see \cite{Gaffney}, \cite{Nuno-Ballesteros-Tomazella}. In particular,  one has:
         $\mu(C,o)=m_1(C,o)-mult(C,o)+1,$ \cite[Theorem 2.6]{Nuno-Ballesteros-Tomazella}. (See also \cite{Brieskorn-Greuel}, \cite{Le}.)
\bex
For plane curves ($c=1$) this one-point scheme is $V(f,\di_y f)\sset (\C^2_{xy},o).$ Its degree is the classical kappa invariant,
 $\kappa(C,o)\!=\!\mu(C,o)\!+\!mult(C,o)\!-\!1,$  see e.g. \cite[pg.212]{Greuel-Lossen-Shustin}.
\eex

\subsection{Preimage of a thin set is thin (lifting tangent arcs to tangent arcs)}\label{Sec.Preliminaries.Thin.Sets}
 The convenient covering $(\C^N,o)\supset X \stackrel{\pi}{\to}(\C^n,o)$ is a Lipschitz map. But $\pi$ is far from being locally bi-Lipschitz.
 In particular, $\pi$ can increase the tangency order of arcs.
Yet, a weaker property holds.

 Recall, a subanalytic germ $Y$  is called {\em thin} if $dim_\R Y> dim T_Y$, see \S\ref{Sec.Preliminaries.Tangent.Cone}.
  The following statement is well known.
\bel\label{Thm.Preimage.of.Thin.Set} Let $\pi:X \to (\C^n,o)$ be a convenient covering and
 $Y\sseteq X  $  a subanalytic sub-germ with connected link.
\bee
\item The \!tangent \!cones \!satisfy: $dim_\R T_Y\!=\!dim_\R  T_{\pi(Y )}.$
In particular, $Y$ is thin iff $(\pi(Y),o)$ is thin.
\item Suppose $T_{(\pi(Y),o)}$ is a (real) half-line.
 \bee
 \item Then the tangent cone $T_Y\sset \C^N$ is a half-line,
  i.e. any two arcs on $Y$ are {\rm outer}-tangent.
  \item Moreover, any two arcs on $Y$ are {\rm inner}-tangent, i.e.
   $tord_{Y}(\ga_1,\ga_2)>1.$
\eee
\eee
\eel

\beR
 Suppose the arcs $\ga_1,\ga_2\sset Y\sset(\R^N,o)$ are real-analytic, i.e. their length-parametrization is by power series with only integer exponents.
  Then their outer tangency order (in $\R^N$) is an integer.
One might hope that the inner tangency order will be an integer and Part 2.b will give  a stronger bound: $tord_{Y }(\ga_1,\ga_2)      \ge2.$
 This does not hold. Take  a real-analytic germ  $Y\sset (\R^3,o)$ that is not LNE, and such that $T_Y$ is a half-line.
  Suppose the (real-analytic) arcs $\ga_1,\ga_2\sset Y$ realize the non-LNE assumption, i.e. $tord_Y(\ga_1,\ga_2)<tord_{\R^3}(\ga_1,\ga_2)= 2.$
   Then $1\le tord_{Y}(\ga_1,\ga_2)<2.$

   An explicit example is the set  $Y\sset \R^3$ defined by the equation $(x^2+y^2)^2=z^3(x^2-y^2)+\ep^2z^7,$ with $0<\ep\ll1.$
  (It is instructive to draw its section by $z=1$.)  The tangent cone is the $\hz$-axis. The two arcs of the intersection $Y\cap \{x=0,z\ge0\}$ are real-analytic. Their parametrization is $(0,\frac{\pm 2\ep z^2}{\sqrt{1+\sqrt{1+z}}},z).$
   Their  outer contact is 2.  But their inner contact is $\frac{3}{2}.$
\eeR

\bcor\label{Thm.Lifting.Tangent.Arcs}
 Suppose two     arcs   $\ga_1,\ga_2\sset (\C^n,o)$ are tangent. For any preimage $\ga^X_1\sset X\sset (\C^N,o)$ there exists a  preimage $\ga^X_2\sset X$ satisfying
 $tord_X(\ga^X_1,\ga^X_2)>1.$
 \ecor
\bpr As the arcs are tangent, we can assume the parameterizations $\ga_i(t)=(t,\uh_i(t)),$ where all the entries of $\uh_i(t)$ are $o(t).$
  Build the H\"older triangle on these arcs,
  \beq\label{Eq.Hoelder.triangle}
  H\ddot{o}ld[\ga_1,\ga_2]:=\cup_{s\in [0,1]} \big[ s\ga_1+(1-s)\ga_2\big]:=
  \cup_{s\in [0,1]} \big[t,s\cdot \uh_1(t)+(1-s)\cdot \uh_2(t) \big]\sset \C^n.
  \eeq
By its construction $H\ddot{o}ld[\ga_1,\ga_2]\sset \C^n$ is a subanalytic subset (with connected link) whose tangent cone is the (real) half-line.
 For a fixed lifting $\ga^X_1\sset X$ we lift the curve $Link[H\ddot{o}ld[\ga_1,\ga_2]]$ (homeomorphically) to
 a curve $Link^X[H\ddot{o}ld[\ga_1,\ga_2]]\sset X$ that starts at $Link[\ga^X_1].$
  This defines the lifting  $H\ddot{o}ld^X[\ga_1,\ga_2]\sset X.$
  Now apply lemma \ref{Thm.Preimage.of.Thin.Set}, with $Y= H\ddot{o}ld^X[\ga_1,\ga_2]$.
 \epr

\beR This corollary fails if any of its assumptions are omitted.
\bee[\bf i.]
\item (Non-convenient covering) Take the covering $\C^3_{xyz}\supset X:=V(z^p-x)\to \C^2_{xy}$ and the parameterized arcs
 $\ga_1(t)=(0,t),$ $\ga_2(t)=(t^l,t)$ for $1<l\le p.$ Then the only liftings are $\ga^X_1(t)=(0,t,0)$ and
  $\ga^X_2(t)=(t^l,t,\om\cdot t^{\frac{l}{p}}),$ where $\om^p=1.$ And $1=tord_{\C^3}(\ga^X_1,\ga^X_2)$
  Thus
  $1=tord_X(\ga^X_1,\ga^X_2)< tord_{\C^2}(\ga_1,\ga_2)=l.$

\item (Non-finite map) Let $\C^3_{xyz}\supset V(zy)\to \C^2_{xy}.$ Take $\ga_1(t)=(t,o)$ and $\ga_2(t)=(t,t^2),$ thus $tord_{\C^2}(\ga_1,\ga_2)=2.$
 For the lifting $\ga^X_1(t)=(t,0,t)$ there does not exist $\ga^X_2$ satisfying: $tord_{\C^3}(\ga^X_1,\ga^X_2)>1.$

\item The conclusion   cannot be strengthened to $tord_X(\ga^X_1,\ga^X_2)=tord_{\C^n}(\ga_1,\ga_2)$. For example, let $f(x,y,z)=z^2+x(x^p+y^p)$, $p>1$.
 The projection $\C^3\supset V(f)\to \C^2_{xy}\supset \De$ is ramified over $\De=V(x(x^p+y^p))$.
 Take two arcs, $\ga_1=(0,t)\sset \De$ and $\ga_2=(t^q,t)$, with $q>p$.
  Then
  \[
  tord_{\C^n}(\ga_1,\ga_2)=q,\quad \ga^X_1=(0,t,0),\quad \ga^X_2=(t^q,t,\pm\sqrt{t^{q+p}(1+t^{q-p})}).
  \]
  Therefore $tord_X(\ga^X_1,\ga^X_2)\le \frac{q+p}{2}<q$.
  \eee
\eeR

\subsection{Fast cycles}\label{Sec.Preliminaries.Fast.Cycles}
  A fast cycle is (roughly speaking) a thin germ that does not admit a hornic deformation-retraction onto a non-thin germ.
\bed\label{Def.Fast.Cycles} Fix a (subanalytic) subgerm $Y\sset X.$
\bee
\item
A subanalytic neighborhood $Y\!\sseteq \!\cU(Y)\!\sseteq \!X$ is  \underline{hornic}  if
 $dif\!f.ord_X(Y,\cU(Y))\!>\!1,$ see \eqref{Eq.Order.of.difference}.

\item $Y$ is called \underline{linkwise retractible} to a subgerm $Y'\sset Y$ if there exists a deformation-retraction $Y\rightsquigarrow Y'$ (inside $X$) that retracts the links, i.e.   $ Link_t[Y]\rightsquigarrow Link_t[Y']$ for all $0<t\ll1.$
\item A thin subgerm $Y\sset X$ is called a \underline{fast cycle} if    $Y$ does not admit
 a hornic neighborhood, $Y\sseteq \cU(Y)\sset X,$ that linkwise retracts to a subgerm $Y'\sset \cU(Y)$ with $dim_\R(Y')=dim_\R(T_{Y'})\le dim_\R T_Y.$

\item Let $Y\sset X$ be a fast cycle. Suppose $Y$ admits a hornic linkwise-retraction to $Y'\sset X,$ with $dim_\R Y'=1+l,$ and this $l$ is the minimal possible.
     Then $Y$ is called \underline{a fast cycle of dimension $l$}.

A fast cycle whose link is $S^1$ is called {\rm a fast loop}.
\item The exponent of a fast cycle is $dif\!f.ord_X(Y,T_Y), $ see \!\eqref{Eq.Order.of.difference}.
  \eee
\eed
A remark: these properties depend on the (subanalytic) inner-Lipschitz type of $X$ only, not on a particular embedding $X\hookrightarrow (\R^N,o).$

\bex
\bee[\bf i.]
\item If a (subanalytic) neighborhood, $Y\sset \cU(Y),$ is hornic, then $T_Y=T_{(\cU(Y),o)}.$ In particular, if $dim_\R\cU(Y)> dim_\R Y,$ then $\cU(Y)$ is a thin set.

 \item A simple  fast loop is the horn $Y=\{x^2+y^2=z^{2\be},\ z\ge0\}\sset \R^3,$ with $\be>1,$ considered inside $X=Ball_1(o)\smin \{(0,0,t)|\ t\ge0\}.$ Note that $Link[Y]$ is contractible inside $Link[X],$ but not hornically.

A simple fast cycle with $dim_\R T_Y>1$ is the cylinder over the horn,  $Y=\{x^2_1+x^2_2=x^{2\be}_3,\ x_3\ge0\}\sset \R^n,$
 considered inside the cylinder $X=\R^{n-3}\times [Ball_1(o)\smin \{x_1=0=x_2,\ x_3\ge0\}].$

\item Suppose $X $ is a hornic neighborhood of its arc $\ga,$ and moreover, $X $ linkwise retracts to $\ga.$ Then $X$ has no fast cycles.
\\For non-hornic neighborhoods $X $ can linkwise retract to $\ga$ and yet can have fast cycles. See  example ii.

\item If $\pi_1(Link[X])\neq \{1\}$ then one can expect fast loops on $X.$ However, fast loops/cycles exist also for exotic spheres, with
  $Link[X]\cong S^{2n-1},$ see Example 4.8 of \cite{Kerner-Mendes.Weighted.Homogen}.
\eee
\eex
\bel\cite[\S2]{Kerner-Mendes.Weighted.Homogen}\label{Thm.IMC.germs.have.no.fast.cycles}
An IMC-germ (subanalytic, with closed link) has no fast cycles.
\eel

\subsubsection{}\label{Sec.Preliminaries.Fast.Cycles.vs.IMC}
 For normal complex-analytic  surface germs the converse of Lemma \ref{Thm.IMC.germs.have.no.fast.cycles} holds:
\\{\bf Theorem 7.5 and Corollary 1.8} of \cite{Birbrair.Neumann.Pichon.14}: {\em If $X$ does not contain fast loops, then $X$ is IMC.}

\

For non-normal surface germs the fast loops are not the only IMC-obstructions.\vspace{-0.2cm}
\bex\label{Ex.Fast.Cycles.not.the.only.IMC.obstructions}
Let $X =V(xz(x-y^2))\sset (\C^3,o).$ It has no fast loops. But $X $ is not IMC.
\bpr
The locus $Sing(X)$ consists of 2 lines and a parabola. At each point of $Sing(X)$ the set $X$ is not a topological manifold.
 If $X$ is inner-Lipschitz equivalent to $Cone[Link[X]]\sset (\R^N,o),$ then $Sing(X)$ must be sent to $Sing[Cone[Link[X]]]=Cone[Sing[Link[X]]].$
  But the parabola  $V(z,x-y^2)\sset Sing(X)$ is inner tangent to the line $V(z,x).$ Hence their images in $Cone[Sing[Link[X]]]$ must be tangent as well. This gives the contradiction.\epr
\eex

\section{Detecting fast loops via the discriminant of covering}\label{Sec.Coverings.Discriminant.Criteria}

\subsection{The covering data}\label{Sec.Covering.Criteria.Ramification.Data}
Take  a reduced complex-analytic   germ  $X\! \sset \!(\C^{n+c}\!,\!o),$  possibly with
\\\parbox{11.5cm}{non-isolated singularity, not necessarily a complete intersection.
Take a convenient covering,  of degree $p.$  Take the discriminant $\De$ with its reduced structure.
 Below $n\!=\!2,$ and we assume: $\De\sset (\C^2,o)$ is a curve germ, see \S\ref{Sec.Preliminaries.Convenient.Coverings}.

\ } \quad\quad
  $\bM (\C^{n+c},o)\supset X  \supset Crit(\pi)
  \\ \quad\quad\searrow \quad \swarrow\pi \quad\swarrow\pi_| \quad
  \\\quad \quad   (\C^n,o) \supset\De\ \ \underbrace{=\cup \De_k }_{for\ n=2}\eM$
\bei
\item Take the tangential decomposition $\De=\cup \De_k.$ Here
 $T_{\De_k}$ is one (possibly multiple) complex line, and $T_{\De_k}\cap T_{\De_i}=o$ for $i\neq k$.
 (But each $\De_k$ can be further reducible.)
\\
Usually we assume that $\De $ is not an ordinary multiple point. Thus $mult(\De_k)\ge2$  for some $k.$

\item  For each $k$ with $mult(\De_k)\ge2$ we take the irreducible decomposition into (complex) branches, $\De_k=\cup  \De_{k,i}.$ They  are all tangent,
   their common tangent (complex) line  is $T_{\De_k}.$

 For each $\De_{k,i}$ let $q_i$ be the total ramification index of the covering  $\pi: X \to  (\C^2,o)$
   over $\De_{k,i}.$ Namely, $q_i=p-\sharp|\pi^{-1}(x)|$ for $x\in \De_{k,i}\smin o.$
   
For generic projections, and $X$ with an isolated singularity, the ramification is minimal, $\{q_i=1\}.$ But we often take convenient non-generic projections.

\item Blowup the origin,  take the strict transforms and exceptional divisors, see the diagram.
\\\parbox{9.1cm}{Note that the curve germs $\{\tDe_j\}$ are now disjoint.

Fix some $\De_k,$ its tangent cone   $T_{\De_k}$ is  just one line. Thus the intersection $\tDe_k\cap \P^1 $ is just one point.  Denote it by $ o_k.$  The fibre $\tpi^{-1}( o_k)\sset \tX$ is a finite set of points, of cardinality
 $\sharp\tpi^{-1}( o_k)\le  p.$
}\hspace{0.3cm}
$\bM\hspace{0.6cm} \P^{1+c}\sset Bl_o(\C^{2+k},o)\supset\tX\supset \tpi^{-1}(o_k)
\\\hspace{3cm}\rotatebox{20}{$\longleftarrow$}\tpi \hspace{1.2cm} \searrow
\\\quad\quad \P^1\sset Bl_o(\C^2,o)\supset\tDe=\coprod \tDe_j\ni o_k\eM$
\\Some of them can be not on a critical locus.
\item
Accordingly $(\tX,\tpi^{-1}( o_k))$ is a (reduced) multi-germ. The number of its connected components is just $\sharp\tpi^{-1}( o_k).$
  Each connected component of $(\tX,\tpi^{-1}( o_k))$ can be further reducible (even if $X $ is irreducible), see Example \ref{Rem.After.Obstruction.Thm}.
      Denote the total number of the connected  components of the multi-germ $(\tX,\tpi^{-1}( o_k))\smin \P T_{X }$ by $r_k.$
       Thus $1\le \sharp\tpi^{-1}( o_k)\le r_k\le p.$

In Step 2 of the proof of theorem \ref{Thm.Fast.Loop.criterion.discriminant.surfaces} we identify $r_k$ also as the number of conected components of $\pi^{-1}Ball,$ where $Ball\sset \C^2$ is a small section transversal to some $\De_k.$
\eei
  \bex\label{Rem.After.Obstruction.Thm}
\bee[\bf i.]
\item If $X $ has an isolated singularity, then $\tX\smin \P T_{X }$ is smooth. Then $r_k$ coincides with the total number of
      irreducible components of the multi-germ $(\tX,\tpi^{-1}( o_k))\smin \P T_{X }.$

\item Cases with $r_k>\sharp\tpi^{-1}( o_k)$ occur often, even if $X $ is irreducible, with an isolated singularity. As the simplest case consider the surface germ $X=V(z^2-f(x,y))\sset (\C^3,o).$ Suppose $ord(f)=2d\ge4$
 and moreover $f(x,y)$ contains the monomial $x^{2d}.$ Blowup the origin, the space $Bl_o(\C^3)\sset \C^3\times\P^2$ is defined by the condition
  $ (x,y,z)\sim[\si_x:\si_y:\si_z].$ Take the chart with $\si_x=1$, its coordinates are $x,\si_y,\si_z.$ (Thus  $y=x\si_y,$ $z=x\si_z.$)
   The strict transform $\tX$ is then defined by $\si^2_z+x^{2d-2}(1+g(x,\si_y))$, where $g(x,\si_y)\in (x,\si_y).$
    This surface germ is reducible, hence $r_k=2>\sharp\tpi^{-1}( o_k)=1.$
    \eee
\eex

\subsection{Fast loop criterion for surface germs}\label{Sec.Covering.Criteria.Surface.Germs}
 (Keeping the notations of \S\ref{Sec.Covering.Criteria.Ramification.Data}.)

\bthe\label{Thm.Fast.Loop.criterion.discriminant.surfaces}
 The germ  $X $ has no fast loops iff one the following holds:
\bei
\item either $\De$ is one point or a curve germ that is an ordinary multiple point (i.e. all its branches are smooth and pairwise non-tangent);
\item or $\De$ is a curve germ and
 each non-smooth tangential component  $\De_k\sset \De$ (with   $mult(\De_k)\ge2 $  and the  data  $\De_k=\cup_i \De_{k,i}$, $q_i$, $r_k$, as above) satisfies:
  $p > r_k-1+\sum_i q_i\cdot  mult(\De_{k,i}).$
\eei
\ethe
\bpr We realize $Link[X]$ as $X\cap V(|x_1|^2+|x_2|^2-\ep), \ 0<\ep \ll 1$, by the convenience of the covering.

 A fast loop (if any) must lie in a hornic neighborhood of the critical locus, $Crit(\pi)\sset X.$
  Therefore in the case ``$\De$ is one point" (and hence $Crit(\pi)$ is a finite set) there can be no fast loops.

 Below we assume that the critical locus is a curve germ, and hence also $\De.$
 Then the image of a fast loop must lie in a hornic neighborhood  $\cU(\De_k)$ for some tangential component $\De_k\sset \De.$ Moreover, this image cannot lie in $\cU(\De_k)$ for a smooth component $\De_k.$ Indeed, if $mult(\De_k)=1,$ then each component of $X$ lying over  $\cU(\De_k)$ is
  linkwise-contractible, see definition \ref{Def.Fast.Cycles}.
 Therefore below we consider only tangential components with $ord\ \De_k\ge2.$
 \bei
\item
In Step 1 we show: $Link[X]$ contains a non-contractible loop ``over $\cU(\De_k)$" iff $X$ has a fast loop in a hornic neighborhood of $\pi^{-1}\De_k.$
 Here the direction $\Lleftarrow$ is immediate.

\item In Step 2 we translate the condition ``$Link[X]$ contains a non-contractible loop over $\cU(\De_k)$"
 into ``$p\le r_k-1+\sum_i q_i\cdot  mult(\De_{k,i})$".
\eei

\bee[\bf Step 1.]
\item

Take the $\hx_1$-axis as the tangent of $\De_k\sset \C^2_{x_1x_2}.$
We work hornically near $\De_k$.
 For each $0<t_o\ll 1$ take
 the section $V(x_1-t_o)\cap \De_k$. These are several points, their number is $mult(\De_k)$.

 As the parameter $ t\in \R_{>0}$ goes to zero  these points draw the subanalytic arcs, $\{\ga_i\}_i\sset \De_k\sset \C^2$.
   Each of these arcs has the length-parametrization $(t,g_i(t))$,
  here $g_i(t)=o(t)$ because $T_{\De_k}=Span_\C(\hx_1)$. In particular,   these arcs $\{\ga_i\}$ are all tangent in $(\C^2,o).$
   For each  pair of these arcs we take the H\"older triangle, as in \eqref{Eq.Hoelder.triangle},
  \beq
  H\ddot{o}ld[\ga_i,\ga_j]:=\cup_{s\in [0,1]} \big[ s\cdot \ga_i+(1-s)\cdot \ga_j\big]:=
  \cup_{s\in [0,1]} \big[t,s\cdot g_i(t)+(1-s)\cdot g_j(t)\big]\sset \C^2.
  \eeq
  This defines the tangent foliation on each $  H\ddot{o}ld[\ga_i,\ga_j].$
  Indeed, any two arcs   $\ga_s,\ga_{s'}\sset  H\ddot{o}ld[\ga_i,\ga_j]$ intersect at $o$ only and are tangent.
 (Because $g_i(t)=g_j(t)$ occurs when either $i=j$ or $t=0$.)

   The union $\C^2\supset H\ddot{o}ld:=\cup_{i\neq j} H\ddot{o}ld[\ga_i,\ga_j]$
 is the ``H\"older polyhedron". This is a subanalytic set and its tangent cone is the real line $Span_\R(\hx_1).$
 Its link is connected and can be realized as $H\ddot{o}ld\cap V(x_1-t_o).$

\quad The lifted H\"older polyhedron is the full pre-image:
\beq
H\ddot{o}ld^X:=\pi^{-1}H\ddot{o}ld=\cup_{i\neq j}\pi^{-1}( H\ddot{o}ld[\ga_i,\ga_j])\sset X.
\eeq
Puncturing at $o$ we get  several connected components, $H\ddot{o}ld^X\smin o=\coprod_l (H\ddot{o}ld^X_l\smin o).$
   Each $H\ddot{o}ld^X_l$ is a (real) subanalytic surface, foliated by the preimages of the foliations on $H\ddot{o}ld[\ga_i,\ga_j]$ for all $i,j.$
 And $Link[H\ddot{o}ld^X_l]$ is connected.
  Lemma \ref{Thm.Preimage.of.Thin.Set} gives: the tangent cone of each $ H\ddot{o}ld^X_l$ is  a real line, and on each $ H\ddot{o}ld^X_l$  the arcs   of our foliation are inner-tangent.

\

  Inside the complex line $V(x_1-t_o)\sset \C^2$
 take an open ball $Ball$ satisfying:
\beq
\De_k\cap V(x_1-t_o)\sset Ball \quad \text{ and } \quad (\De\smin \De_k)\cap Ball=\empty.
\eeq
 Such a ball exists  for $0<t_o\ll1$ as   every branch of $\De\smin \De_k$ is non-tangent to $\De_k.$

   Suppose the full preimage $\pi^{-1}Ball\sset X\cap V(x_1-t_o)$ contains a (compact, connected) non-contractible loop $Z$.
      By deforming/shrinking $Z\sset \pi^{-1}Ball $ we can assume $\pi(Z)\sset H\ddot{o}ld\cap Ball.$
      Indeed, fix a simply-connected path  $\eta\sset H\ddot{o}ld\cap Ball$ that visits all the points of  $\De_k\cap V(x_1-t_o).$
      The retraction $Ball\rightsquigarrow\eta$ lifts to $\pi^{-1}Ball\rightsquigarrow \pi^{-1}\eta.$

    Therefore $Z$ lies inside  $   H\ddot{o}ld^X \cap V(x_1-t_o).$ Then it lies inside a connected component,
    $Z\sset H\ddot{o}ld^X_l \cap V(x_1-t_o)$ for some $l.$ Finally, one can assume: $Z$ is subanalytic, \S\ref{Sec.Preliminaries.Dictionary}.v.

      Using the (inner-tangent) foliation on that component $H\ddot{o}ld^X_l$ we get        the subanalytic surface germ
      \beq
      \R_{>0}Z:=\cup_{z\in Z}\ga_z\sset H\ddot{o}ld^X_l\sset X.
      \eeq
       Its tangent cone is a real half-line, it lies inside the tangent cone of $ H\ddot{o}ld^X_l.$
        Thus         $\R_{>0}Z$ is a thin germ.
        By our construction the germ  $\R_{>0}Z$  is not linkwise-retractible inside $\pi^{-1} Cone[Ball],$ see \S\ref{Sec.Preliminaries.Fast.Cycles}. Hence it is not linkwise retractible inside any hornic neighborhood of $T_{\R_{>0}Z}$
 Therefore $\R_{>0}Z$ is the promised fast loop.

 \

\item  It remains to check whether a non-contractible (connected, compact) loop $Z\sset \pi^{-1}Ball$ exists. Namely (in terms of the first homology), whether  $h^1(\pi^{-1}Ball)>0.$
 Therefore it remains to show the equivalence:
 \beq\label{Eq.inside.proof}
 h^1(\pi^{-1}Ball)>0 \quad\quad\quad  \text{ if and only if  }\quad\quad\quad   p\le r_k-1+\sum_i q_i\cdot  mult(\De_{k,i}).
 \eeq
  Thus we compute  the Euler characteristic  $\chi(\pi^{-1}Ball)$ and  the number of connected components
     $h^0(\pi^{-1}Ball).$

\bei
\item
The map $X\supset \pi^{-1}Ball\stackrel{p:1}{\to}   Ball\sset V(x_1-t_o)\sset \C^2$ is a covering of complex  analytic  curves.
 If $X$ has a non-isolated singularity then the complex-analytic curve $\pi^{-1}(Ball)$ can be singular.
 The covering is ramified over
 the discriminantal points $\De_k\cap V(x_1-t_o) $.
 For the branch decomposition $\De_k=\cup \De_{k,i}$ we get  $\sharp|\De_{k,i}\cap V(x_1-t_o)|=mult(\De_{k,i})$, and these points are of ramification index $q_i$.
 Riemann-Hurwitz  formula applied to this covering gives:
 \beq
 \chi(\pi^{-1}Ball)=p \cdot \chi(Ball)-\sum q_j=p -\sum q_i\cdot  mult(\De_{k,i}).
\eeq

 \item
We compute the number of connected components:  $h^0(\pi^{-1}Ball)\!=\!r_k.$ Blowup the
\\
\parbox{8cm}{origin to get the diagram,  as in  \S\ref{Sec.Covering.Criteria.Ramification.Data}.   Take the intersection point
    $o_k\!:=\tDe\!\cap \!\P^1\!\cap \!\widetilde{Cone[Ball]}\!=\!\tDe_k\!\cap \!\P^1.$
  Take the pre-images $\tpi^{-1}(o_k)\sset \P T_X\sset \tX.$
}\quad \setlength\arraycolsep{1pt}
$\bM   \coprod\tDe_k\sset Bl(\C^2,o)&\stackrel{\tpi}{\leftarrow}\tX  \sset& Bl_o(\C^N,o)
 \\ \downarrow\hspace{1.3cm}\downarrow\si&\downarrow &\downarrow
 \\\cup\De_k\sset (\C^2,o)&\stackrel{\pi}{\leftarrow} X\sset&(\C^N,o)
 \eM$

  We get the multi-germ $ (\tX,\tpi^{-1}(o_k)).$
 By definition, $r_k$ is the number of connected components of
  $(\tX,\tpi^{-1}(o_k))\smin \P T_{X }=:\coprod \tY_l.$ Here each $\tY_l$ is a connected (locally-closed) complex-analytic surface.

 Starting from $Ball\sset V(x_1-t_o)\sset \C^2,$  we get $(\si\circ\tpi)^{-1}Ball\sset \widetilde{V(x_1-t_o)}\cap\tX.$
  As we work with the multi-germ $(\tX,\tpi^{-1}(o_k))$  we can assume:
   $(\si\circ\tpi)[\widetilde{V(x_1-t_o)}\cap (\coprod \tY_l)]\sset Ball.$
  Therefore we can ignore all the components of $\De\smin \De_k,$ i.e. we can assume $\De=\De_k.$

\

  It is enough to prove: $\tY_l\cap ((\si\circ\tpi)^{-1}Ball)$ is connected for each $l.$
    This will imply: the number of connected components of $\coprod \tY_l \cap ((\si\circ\tpi)^{-1}Ball)$ coincides with the number of connected components of $\coprod \tY_l.$ And therefore:
    \beq
    h^0(\pi^{-1}Ball)=  h^0((\si\circ\tpi)^{-1}Ball)=h^0((\coprod \tY_l) \cap ((\si\circ\tpi)^{-1}Ball))=h^0(\coprod \tY_l)=r_k.
  \eeq
 Fix two points, $\tp_1,\tp_2\in \tY_l\cap (\si\circ\tpi)^{-1}Ball.$ We can assume (w.l.o.g.) $\tp_1,\tp_2\not\in Crit(\tpi).$
  Take a path $[\tp_1\stackrel{\teta}{\rightsquigarrow}\tp_2]\sset \tY_l.$ We want to deformation-retract $\teta$ into a path inside
  $\tY_l\cap (\si\circ\tpi)^{-1}Ball.$ Take the image $\eta:=(\si\circ\tpi)(\teta)\sset \C^2\smin \{o\}.$ It connects the points
   $\si(\tpi(\tp_1)),\si(\tpi(\tp_2))\in Ball.$
  \bee[\bf i.]
  \item Suppose $\teta\cap Crit(\tpi)=\empty.$ Then we can assume $\eta\cap \De=\empty,$ i.e. $\eta\sset \C^2\smin \De.$ Then $\eta$ can be deformation-retracted   (inside $\C^2\smin \De$, and preserving the end-points)  to a path inside $Ball\smin (Ball\cap \De).$

      This deformation lifts uniquely to the deformation of $\teta$   inside $\tX\smin Crit(\tpi)$, while preserving the points $\tp_1,\tp_2.$  The deformed path,
         $[\tp_1\stackrel{\teta_1}{\rightsquigarrow}\tp_2]$ lies fully inside $\tY_l\cap (\si\circ\tpi)^{-1}Ball.$
\item The path $\teta$ can happen to intersect the critical locus $Crit(\tpi).$ Moreover, $\teta$ can be non-deformable off   $Crit(\tpi).$
 E.g. let $Y_l$  be connected but reducible, with $\tp_1,\tp_2$ lying in distinct components. Then $\teta$ inevitably intersects the locus
  $Sing(\tY_l)\sseteq Crit(\tpi).$

 \quad  We claim: $\teta$ can be  deformed (inside $\tY_l$) to a path $[\tp_1\stackrel{\teta_1}{\rightsquigarrow}\tp_2]$ satisfying:
  $\teta_1\cap Crit(\tpi)\sseteq (\si\circ\pi)^{-1}Ball\cap \tY_l.$ Indeed, we can assume (w.l.o.g) that $\teta_1\cap Crit(\tpi)$ is a finite number of points. Therefore the set $\eta\cap \De$ is finite. Over each such intersection point we can slightly deform $\teta$ to get: the intersection $\eta\cap\De$ is transverse.

  Now we slide-move all the points of $\eta\cap \De$ (along $\De$) into $Ball\cap \De.$ This deformation lifts uniquely to a deformation of $\teta.$  (Because $\De\smin o$ is a smooth curve and $\eta\pitchfork \De.$)

  Finally, having reached the path with $\teta_1\cap Crit(\tpi)\sseteq (\si\circ\pi)^{-1}Ball\cap \tY_l,$ one applies  part i. to all the connected components of $\teta_1\smin (\teta_1\cap Crit(\tpi)).$
  \eee
\eei

\

\noindent We have proved: $h^1(\pi^{-1}Ball)=r_k-p+\sum q_i\cdot  mult(\De_{k,i}).$
 This gives equation \eqref{Eq.inside.proof}
  Hence the statement.
\epr
\eee

Combining this theorem with \S\ref{Sec.Preliminaries.Fast.Cycles.vs.IMC}, we get:

\bcor\label{Thm.IMC.criterion.Corollary}
Suppose $X $ is a normal surface germ and $\De$ is a curve germ.
   Then $X $ is IMC iff   all the non-smooth tangential components $\De_k\sset \De$ (with   $mult(\De_k)\ge2 $  and the  data  $\De_k=\cup_i \De_{k,i}$, $q_i$, $r_k$, as above) satisfy:
  $p > r_k-1+\sum_i q_i\cdot  mult(\De_{k,i}).$
\ecor

\subsubsection{An example: $X  \sset (\C^{2+c},o)$ is weighted-homogeneous, with weights $\om_1\le \om_2\le \cdots\le \om_{2+c} $}
 \bcor\label{Thm.IMC.criterion.semi.weighted.homogen.germs}
Suppose  $\om_1=\om_2$ and  the projection $X \to (\C^2_{x_1x_2},o)$  is convenient. Then $X $ has no fast loops.
If, moreover, $X $ is a normal surface germ then $X $ is IMC.
  \ecor
A remark, if $\om_1<\om_2,$ then $X$ necessarily has a fast loop, \cite[\S4]{Kerner-Mendes.Weighted.Homogen}.
\bpr
 The discriminant $\De \sset (\C^2,o) $ of the projection $X \to(\C^2_{x_1x_2},o) $ is weighted-homogeneous with weights $\om_1,\om_2.$ Hence it is homogeneous, i.e. an ordinary multiple point.
 Now invoke (trivially) theorem \ref{Thm.Fast.Loop.criterion.discriminant.surfaces}.
\epr

\subsubsection{An example: surface germs of multiplicity 2}\label{Sec.Examples.IMC.order=2}
Such a germ is necessarily a  surface in $(\C^3,o)$. Thus we start from $X =V(z^2-a_0(x,y))\sset (\C^3,o),$ where $ord(a_o)\ge2.$
 \bcor
 $X $ has no fast loops iff the curve germ $V(a_o(x,y))\sset (\C^2,o)$ is set-theoretically an ordinary multiple point.
 \ecor
For   $X $ with an isolated singularity this gives the IMC criterion.
\bpr
 The covering $X \to (\C^2,o)$ is totally ramified over the discriminantal curve, $\De:=V(a_o(x,y))\sset (\C^2,o),$  taken set-theoretically.
 If $\De$ is an ordinary multiple point, then there is nothing to check in theorem \ref{Thm.Fast.Loop.criterion.discriminant.surfaces}.
 Otherwise $\De$ contains a non-smooth tangential component, $\De_k,$ with $mult(\De_k)\!\ge\!2.$ Then
  $p\!=\!2\!\le\! mult(\De_k)\cdot (2-1)+(r_k-1).$ Hence $X $ has a fast loop.
  \epr

\bex\label{Ex.IMC.of.mult2} (IMC's of multiplicity 2 and right-modality$\le2$)
\bei
\item
Among the simple types (i.e. ADE's) the only IMC's are: $A_1$ ($x^2+y^2+z^2$) and $D_4$ ($z^2+x^3+y^3$).
 This is well known, see e.g. \cite{Birbrair-Fernandes.Neumann.08}.
\item  The only  non-simple IMC germ of multiplicity 2 and right-modality$\le2$ is  $X_9$ ($z^2+x^4+y^4+a\cdot x^2y^2$),
  see the tables in \cite[pg.24]{AGLV}. The next case,  $z^2+x^5+y^5+a_{3,2}  x^3y^2+a_{2,3}  x^2y^3+a_{3,3}  x^3y^3,$  has modality$\ge3.$
\eei\eex

\subsubsection{An example: germs of type ${V(z^p+za_1(x,y)+a_0(x,y))\!\sset\! (\C^3,o),\ p\!\ge\!3}$}
Here $ord (a_o)\!\ge\! p,$ $ord (a_1)\!\ge\! p-1.$
 By \S\ref{Sec.Preliminaries.Crit.And.Discr} we get a simple expression for the discriminant,
$\De=(\frac{a_0(x,y)}{p-1})^{p-1}-(\frac{a_1(x,y)}{p})^p\in \C\{x,y\}.$ This can be not square-free.
 Denote  its reduced version by $\De_{red}.$ Take the lowest order terms, $l.t.(\De),$ resp. $l.t.(\De_{red}).$  These homogeneous polynomials  define  the   tangent-cone
  of $V(\De),$ resp. of $V(\De_{red}).$ Blowup  the origin, $\P^2\sset Bl(\C^3,o) \to (\C^3,o),$ and take the strict transform $\tX\to X.$
 \bcor
 \bee
 \item Suppose $a_1=0.$ Then $X $ has no fast loops iff $V(a_o)\sset (\C^2,o)$ is set-theoretically an ordinary multiple point.
  \item Suppose $a_1\neq0.$ Factorize the lowest order terms, $l.t.(\De)=\prod l_k(x,y)^{m_k},$ resp. $l.t.(\De_{red})=\prod l_k(x,y)^{m_{k,red}}.$
  Then $X $ has no fast loops iff for each $k$ with $  m_{k,red}\ge 2$ the following conditions hold:
  \bee[\bf i.]
  \item       $ m_k\le p.$
  \item The total number of   irreducible components of the analytic germ $\tX$ at all the points of
   $\tX\cap \P^2\cap \widetilde{V(l_k)}$ is at most $p-m_k.$
  \eee
 \eee
 \ecor
A remark: in most cases already the condition 2.i. is violated, i.e. one has $m_k> p.$
\bpr In both cases we should verify the bound  $p> r_k-1+\sum_i q_i\cdot  mult(\De_{k,i})$ of theorem \ref{Thm.Fast.Loop.criterion.discriminant.surfaces}.
\bee
\item  If $a_1=0$ then the covering $X \to (\C^2,o)$ is totally ramified over the curve $V(a_o(x,y))\sset (\C^2,o).$
 The total ramification index is then  $q_i=p-1,$ and $r_k\le p.$
  Then the bound holds iff   $mult(\De_k)=1$ for all $k.$
\item  In this case $m_{k,red}=mult(\De_k)$ and
  $\sum_i q_i\cdot  mult(\De_{k,i})=m_k.$ Thus the necessary and sufficient condition to be IMC, i.e. $p> r_k-1+m_k,$ means exactly the conditions i. and ii.
  \epr
\eee
Already in this (very particular) case we get a vast amount of IMC's.
\bex Let $X=V(z^p+za_1(x,y)+a_0(x,y))\sset (\C^3,o),$ with
  $p\ge 3,$ and an isolated singularity. Suppose the power series $a_0(x,y),a_1(x,y)$ are co-prime.
\bei
\item Suppose $\frac{ord(a_1(x,y))}{p-1}>\frac{ord (a_0(x,y))}{p}  .$ Then $l.t.(\De)=l.t.(a_0(x,y)^{p-1}).$ Thus
$X $ is IMC iff the (set-theoretic) curve germ $V(a_0(x,y))\sset (\C^2,o)$ is an ordinary multiple point.
\item Suppose $\frac{ord(a_1(x,y))}{p-1}<\frac{ord (a_0(x,y))}{p}  .$ Then $l.t.(\De)=l.t.(a_1(x,y)^{p}).$
 Thus $X $ is IMC iff the (set-theoretic) curve germ $V(a_1(x,y))\sset (\C^2,o)$ is an ordinary multiple point.
\eei
\eex

\subsection{Fast loops imposed by singular points of the tangent cone}

  Suppose a surface germ  $X\sset (\C^{2+c},o)$
   is a strictly complete intersection with {\em isolated} singularity, \S\ref{Sec.Preliminaries.Tangent.Cone.C.I.}.
       The  projectivized tangent cone  is then a (globally) complete intersection curve, $\P T_X\sset \P^{1+c}\sset Bl_o(\C^{2+c},o).$
       It can be reducible, non-reduced, i.e. can have multiple components.
        Occasionally we take the underlying set-germ (with reduced structure), $\P T^{red}_X.$
        One has the total degree,   $deg[\P T_{X }]=p,$ and the set-theoretic degree,
        $deg[\P T^{red}_{X }]\le p.$

Call a point $pt\in \P T_X$ {\em reduced} if the curve germ $(\P T_X,pt)$ is reduced.
 Thus $(\P T_X,pt)$ is an ICIS, and its Milnor number, $\mu(\P T_X,pt),$ is well defined.

Call a hyperplane $\P^c\sset \P^{1+c}\sset Bl_o(\C^{2+c},o)$ {\em non-tangent} to the curve $\P T_X$ if $\P^c$ is transverse to the tangent cone of the germ $(\P T_X,pt)$ at  each intersection point $pt\in \P^c\cap \P T_X.$
\bex
\bee[\bf i.]
\item The generic hyperplane $\P^c\sset \P^{1+c} $ is non-tangent.  It intersects $\P T^{red}_X$ transversally and only at smooth points of $\P T^{red}_X.$
\item Suppose the germ $(\P T_X,pt)$  is reduced and singular.
 Among all the hyperplanes  through $pt,$ take the generic one. It is non-tangent to $\P T_X.$
\item
 More generally, take some reduced singular points, $pt_1,\dots,pt_j\in Sing(\P T_X),$ with $j\le c.$ The linear system of hyperplanes through these points has dimension$\ge c+1-j.$
 And the generic member of this family is non-tangent to $\P T_X$ at the points $pt_1,\dots,pt_j,$ unless these points are in degenerate position with respect to the curve $\P T_X.$
  \eee\eex

   \bthe\label{Thm.Fast.Loops.via.Tangent.Cone.surfaces}
   \bee
   \item If   $X$ is   IMC and the curve $\P T_X$ is reduced then $\P T_X$ is  smooth. (I.e. $T_X$ has an isolated singularity, i.e. $X$ is a ``generalized  ordinary multiple point".)

   \item 
   Suppose a  hyperplane $\P^c\sset \P^{1+c}$ is non-tangent to $\P T_X,$ and contains reduced singular points  $\{pt_i\}$ of the curve $\P T_X.$
 If their Milnor numbers satisfy $\sum \mu(\P T_X,pt_i)\ge 1+deg[\P T_{X }]-deg[\P T^{red}_X],$ then $X$ has a fast loop.

 In particular, if $\P T_X=\P T_X^{red}$, any (reduced) singular point of $\P T_X \cap \P^c$ produces a fast loop.

   \eee
\ethe
The converse to part 1 is Example 3.8 in \cite{Kerner-Mendes.Weighted.Homogen}.
\bpr
Part 1 follows straight from Part 2. Part 2 is proved in steps.
\bei
\item
Take the (convenient) projection $\C^{2+c}\to \C^2$ associated to $\P^c.$ Then $\P^c$ is sent to a point $o_k\in \P^1\sset Bl_o(\C^2,o).$
 Take the corresponding tangential component of the discriminant, $\De_k\sseteq\De.$
 To use theorem \ref{Thm.Fast.Loop.criterion.discriminant.surfaces}
we verify the condition  $p \le r_k-1+ \sum q_i\cdot mult(\De_{k,i}).$
   First we bound $r_k\ge \sharp \tpi^{-1}(o_k).$
 \item
 We express    $\sum q_i\cdot mult(\De_{k,i}) $ as the multiplicity of the polar curve, and use the formula of \S\ref{Sec.Preliminaries.Crit.and.Discr.Polar.Multipl}.
\eei
\bee[\bf Step 1.]
\item 
 The projection $(\C^{2+c},o)\stackrel{\pi}{\to}(\C^2_{x_1x_2},o)$ lifts to
 $Bl_o(\C^{2+c},o)\supset \P^{1+c}\smin \widetilde{V(x_1,x_2)}\to \P^1\sset Bl_o(\C^2,o).$ By a preliminary $GL(\C^{2+c})$ transformation one can assume:
\beq
\P^c=\widetilde{V(x_1)}\cap \P^{1+c},\quad\quad
\widetilde{V(x_1,x_2)}\cap \P T_X=\empty,
\eeq
\[
\hspace{2cm}\text{the intersection }\widetilde{V(a_1x_1+a_2x_2)}\cap \P T_X \text{ is finite for all} [a_1:a_2]\in \P^1.\]
Indeed, once the first condition is satisfied, one uses the remaining freedom $x_i\to x_i+l_i(x_1\dots x_{2+c})$ for $i=2,\dots, 2+c,$ and generic linear forms $\{l_i\}.$

Thus $\pi$ sends $\P^c\smin \widetilde{V(x_1,x_2)}$ to the point $o_k:=[0:1]\in \P^1.$ We get the diagram:

\beq\bM
&&(\C^{2+c},o)&\leftarrow&Bl_o(\C^{2+c},o)&\supset&\P^{1+c}
\\&&\cup   &&   \cup    &&\cup
\\
Crit(\pi)&\sset &  X  &\leftarrow&\tX  &\supset&\P T_X&\supset & \tpi^{-1}(o_k)&\supseteq \{pt_i\}
\\ \downarrow&&\pi\downarrow   && \tpi \downarrow   && \downarrow \tpi_|&& \downarrow
\\\De&\sset&(\C^2,o)&\leftarrow  & Bl_o(\C^2,o)&\supset & \P^1 & \ni &o_k&
\eM
\eeq
Here are  the needed properties of this diagram.

 \bee[\bf i.]
 \item
The projection $\pi $ is convenient, as   $T_X\cap V(x_1,x_2)=(o),$ as $\widetilde{V(x_1,x_2)}\cap \P T_X=\empty.$

\item The projection $\tpi$ is finite, as $\pi$ is convenient. (By part 1 of Lemma \ref{Thm.Diagram.Blowup.Lift.Covering}.)
\item The curve
$\P T_X\sset \P^{1+c}$ is a complete intersection of the (total) degree $p.$
    Therefore, $\tpi_|$ is a  (flat, ramified) $p:1$ covering. (Part 2 of Lemma  \ref{Thm.Diagram.Blowup.Lift.Covering}.)

\item
The fibre $\tpi^{-1}(o_k)\sset \P T_X$  contains the reduced singular points $\{pt_i\}.$
 Its cardinality   is bounded, $\sharp \tpi^{-1}(o_k) \le p-\sum_{\tilde{o}_j\in \tpi^{-1}(o_k)} (mult(\P T_{X },\tilde{o}_j)-1).$
 As the hyperplane $\P^c$ is non-tangent to $\P T_X,$ this bound is an equality.

\noindent Altogether, we get:
\beq\label{Eq.proof.of.Thm.fast.loops.from.PTX}
r_k\ge \sharp \tpi^{-1}(o_k) = p-\sum_{\tilde{o}_j\in \tpi^{-1}(o_k)} (mult(\P T_X,\tilde{o}_j)-1)=
\eeq\vspace{-0.3cm}
\[
\hspace{7cm}
=p-\sum (mult(\P T_X,pt_i)-1)-deg[\P T_X]+deg[\P T^{red}_X].\]
 \item The projection $Crit(\pi)\to \De$  lifts to the (finite) projection
 $Crit(\tpi)\to \tDe\cup \P^1.$  Here $\tX$ is locally-complete intersection, thus $Crit(\tpi)$ as well (part 1 of Lemma \ref{Thm.Discriminant.Basic.Properties}).

  We claim: $Crit(\tpi)\!=\!\widetilde{Crit(\pi)}\!\cup\! Non.Red(\P T_X),$  the non-reduced components of $\P T_X.$
  Indeed, $\tpi$ is not a local isomorphism at any point of $Non.Red(\P T_X).$ (Because the local ramification index of $\tpi$ at each point of $Non.Red(\P T_X)$ is   the multiplicity of $Non.Red(\P T_X),$ hence $\ge2.$)
    Vice versa, $Crit(\tpi)$ is a locally complete intersection curve,
   thus $\overline{Crit(\tpi)\smin Non.Red(\P T_X)}$ is of pure dimension 1.
   Hence, besides $\widetilde{Crit(\pi)},$ it can consist only of
   some irreducible components of $\P T_X.$
 But the restriction of $\tpi$   to each irreducible, reduced component of $\P T_X$ is a local isomorphism (over $\P^1$) except at a finite number of points.

\eee

 \item
  Take the image $\tDe_k:=\tpi(Crit(\tpi),\pi^{-1}(o_k))\sset (Bl_o(\C^2,o),o_k).$ This is the strict transform of $\De_k\sset (\C^2,o),$
  a tangential component of $\De.$

   Take the coordinate change $\Phi:x\to x+h(x),$ where the terms $h(x)=(h_1(x),\dots,h_{2+c}(x))$ of order two and generic.
  This lifts to an automorphism $\tPhi\circlearrowright Bl_o(\C^{2+c},o),$ with $\tPhi|_{\P^{1+c}}=Id.$ (Because $\Phi'|_o=Id.$)
   Therefore $\tPhi$ preserves $o_k$, and all the properties of Step 1 are preserved. And applying $\Phi$ the covering
  $X\to (\C^2,o)$ becomes   generic in the following sense: $q_i=1$ for all $i.$

  Therefore $\sum_i q_i\cdot mult(\De_{k,i})=mult(\De_k,o),$
  and we compute this multiplicity.

Pass from $(\C^2,o)$ to $Bl_o(\C^2,o)$ by
  $  mult(\De_k)=  \bl  \tDe_k,\P^1\br.$ By our choice of coordinates, $o_k=[0:1]\in \P^1,$ i.e. the tangent line of $\De_k$ is $Span(\hx_2).$
   Thus at any point of $\tpi^{-1}(o_k)$ the exceptional divisor is defined locally by $x_2=0.$
    Then  $\bl  \tDe_k,\P^1\br=\bl  \tDe_k,\widetilde{V(x_2-t_o)}\br,$ as $\tDe_k\sset Bl_o(\C^2,o)$ is a hypersurface.
    Take the multi-germ $Crit_k(\tpi):=(Crit(\tpi),\tpi^{-1}(o_k)).$
  Then
  \beq
   mult(\De_k)=\bl Crit_k(\tpi),\widetilde{V(x_2-t_o)}\br \ge \bl \overline{Crit_k(\tpi)\smin \P T_{X }  },\P^{1+c}\br\ge
  \sum \bl (Crit(\tpi),pt_i),\P^{1+c}\br.
  \eeq
Both inequalities hold here as $Crit(\tpi)$ (with its scheme structure) is a locally complete intersection (lemma \ref{Thm.Discriminant.Basic.Properties}).

 At the chosen reduced singular  points, $pt_i\in \P T_X\cap \P^c, $
one has:
 $\bl   Crit(\tpi),\P^{1+c}\br_{pt_i}=\mu_i +mult(\P T_X,pt_i)-1.$ (Apply \S\ref{Sec.Preliminaries.Crit.and.Discr.Polar.Multipl} to the  reduced curve germ with $(\P T_{X },pt).$)

Altogether we have: $  mult(\De_k)\!\ge\! \sum (\mu_i+mult(\P T_X,pt_i)-1).$   In particular  $mult(\De_k)\!\ge\!2.$

\eee

\

Combine this with \eqref{Eq.proof.of.Thm.fast.loops.from.PTX} to get:
$ r_k+   mult(\De_k)-p\ge  \sum \mu_i-deg[ \P T_X]+deg [ \P T^{red}_X ].$
 Together with the assumption $ \sum\mu_i\ge 1+ deg[ \P T_X]-deg [ \P T^{red}_X ]$ we get:
  $ r_k + \sum q_i\cdot mult(\De_{k,i})\ge p+1 .$ Hence the statement.
  \epr

\bex\label{Ex.IMC.conditions.imposed.by.PT_X}
\bee[\bf i.]
\item
For  any super-isolated surface singularity in $(\C^3,o)$ the projectivized tangent cone is reduced. In this case $X $  is  IMC iff it is an ordinary multiple point. This was obtained in \cite{Birbrair.Neumann.Pichon.14}, example 15.1.

\item
Suppose a strictly complete intersection surface germ (with isolated singularity) is both LNE and IMC. Then
   $X $ is an ordinary multiple point.   In particular,  $X \cong (T_{X },o).$

    Indeed,  $T_{X }$ of an LNE is reduced, \cite{Fernandes.Sampaio}.
    Thus $\P T_{X }$ is smooth.
 \eee
\eex
\beR
 Part 2 of this theorem does not hold for germs with non-isolated singularities.
  E.g. let $X\sset (\C^3,o)$ be an arrangement of planes (through the origin). Then $X$ is an IMC and LNE, and $\P T_X$ is reduced but non-smooth.
  \eeR

\subsubsection{An example: IMC's   of right modality $\le2$} (Continuing example \ref{Ex.IMC.of.mult2}.)
\bcor\label{Thm.IMC.among.right.modality.le.2}
Among the   surface singularities in $(\C^3,o)$ of right modality $\le2$ the only IMC's are: $A_1,D_4,$ $P_8$ ($x^3+y^3+z^3+a\cdot xyz$), and
  $X_9$ ($z^2+x^4+y^4+a\cdot x^2y^2$).
\ecor
\bpr
The case of corank=2 was treated in  \S\ref{Sec.Examples.IMC.order=2}.
Thus we consider only the cases of corank=3. Their tables  are listed e.g. in   \cite[pg. 24]{AGLV}.
 In particular, the only ordinary multiple points   of modality $\le2$ are $D_4,P_8.$

Suppose $X $ is IMC and  not an  ordinary multiple point. Then the curve $\P T_{X }$ must have multiple components.
 But  all the corank=3 surface singularities of modality $\le2$ have reduced tangent cone.
  (By the direct check of those tables.)
 \epr

\subsection{Detecting fast loops in higher dimensions, $n\ge3$}
 Let  $X\sset (\C^N,o)$ be as in \S\ref{Sec.Preliminaries.Dictionary}.i,  with an arbitrary singularity.
 Take a convenient covering, $(\C^N,o)\supset X\to (\C^n,o)\supset \De.$
 We assume $dim_\C\De=n-1,$ e.g. this holds if $X$ is a complete intersection.

  Fix a smooth surface germ  $\cS\sset (\C^n,o)$, and
 (set-theoretically) restrict the covering,  $X\cap \pi^{-1}\cS\stackrel{\pi|_\cS}{\to} \cS\sset (\C^n,o).$
 Then $X \cap \pi^{-1}\cS$ is a (pure dimensional, possibly reducible) surface, and the projection $\pi|_\cS$ is still a convenient covering (see \S\ref{Sec.Preliminaries.Convenient.Coverings}), of degree $p$.
 Now $\pi|_\cS$ is ramified over
  the curve germ $\De|_\cS:=\De\cap \cS.$
   As in \S\ref{Sec.Covering.Criteria.Ramification.Data} we take a tangential component, $\De|_\cS\rightsquigarrow \De_k=\cup (\De_{k,i},o),$ and the corresponding covering data $r_k,$ $\{q_i\},$ $\{mult(\De_{k,i})\}.$

Observe that the embedded topological type  of the planar curve germ $\De|_\cS\sset \cS\cong (\C^2,o) $   can vary when  the section $ \cS\sset (\C^n,o)$ varies.

\bprop\label{Thm.Fast.Loops.via.discriminant.higher.dim}
   Suppose the section  $\De|_\cS$  has a   tangential component $\De_k$ with
 $mult(\De_k)\ge 2. $   Suppose any deformation $\cS(t)$ of the surface $\cS=\cS(o)$ that preserves the tangent plane
  $T_\cS\sset \C^n $  results in an equisingular family of the planar curve germs,    $\De_k(t)\sset \cS(t).$
   If  $p\le r_k-1+\sum_i q_i\cdot  mult(\De_{k,i}) $  is satisfied for $\cS(0),$ then $X $ has a fast loop.
\eprop
Here $t\in \R^1$ in the family $\cS(t),$ i.e. $t$ is not necessarily small.
\bpr
Restrict  the covering as in the assumption, $X\supset X_\cS\to \cS\supset \De|_\cS.$
 This is still a convenient covering.
 As in the proof of theorem \ref{Thm.Fast.Loop.criterion.discriminant.surfaces} we get a fast loop $\R_{>0}Z\sset  X_\cS.$
  We show that it is a fast loop also inside $X.$

  Any deformation of $\cS$ that preserves the tangent plane $T_\cS$ induces an equisingular family  $\De_k(t).$
   Thus it induces a family of  fast loops, $\R_{>0}Z(t)\sset X_{\cS(t)},$ with no hornic link-wise retraction, \S\ref{Sec.Preliminaries.Fast.Cycles}. Thus $\R_{>0}Z$ could be possibly link-wise retractible  inside $X$ only in a family     with varying tangent planes. But such a retraction is non-hornic.
    Therefore $\R_{>0}Z$ is a fast loop also inside $X$, by \S\ref{Sec.Preliminaries.Fast.Cycles}.
\epr

\beR
The assumption of equisingularity of $\De_k(t)\cap \cS(t)$
is important.
 Without it one could start, e.g. from  a homogeneous $\De\sset (\C^n,o)$, and take an intersection with a tangent smooth surface, to get a singular branch of high multiplicity. This would ensure a germ that is a fast cycle inside $X\cap \pi^{-1}\cS,$ but   admits a hornic linkwise  retraction inside $X.$
\eeR

As a corollary we get: points of high multiplicity of  $\P T_{\De}\sset \P^{n-1}\sset Bl_o(\C^n,o)$  ensure fast loops in $X.$
\bcor\label{Thm.Fast.Loops.via.PT.discriminant.higher.dim}
Let $X\sset (\C^{n+c},o)$ be a complete intersection of dimension $n\ge3$ and multiplicity $p.$ Take a convenient covering $X\to (\C^n,o)\supset \De$ and suppose the total ramification index over each point of $\De$ is at least $q,$ see \S\ref{Sec.Covering.Criteria.Ramification.Data}.
 Suppose a point $o_k\in \P T_\De$  satisfies:
 $mult(\P T_\De,o_k)\ge \frac{p}{q}$ and $mult(\P T_\De,o_k)> mult(\tDe,o_k).$
  Then $X$ has a fast loop.
\ecor
\bpr
 To apply proposition \ref{Thm.Fast.Loops.via.discriminant.higher.dim} we blowup $Bl_o(\C^n,o)\stackrel{\si}{\to}(\C^n,o)$ and
  fix a smooth ($\C$-analytic) surface germ $\cS\sset (\C^n,o)$ satisfying:
\[
 o_k\in \tilde{\cS}\cap \P^{n-1}\sset Bl_o(\C^n,o) \quad  \text{but
$\cS$ is generic otherwise.}
\]
Then $\cS\cap \De$ is a reduced plane curve germ. One of its tangential component is $\De_k:=\si(\tilde\cS\cap \tDe, o_k).$
 We bound its multiplicity, $mult( \De_k)=\bl \tDe_k,\P^{n-1}\br=mult(\P T_\De, o_k)\ge \frac{p}{q}.$ Therefore we have the assumption of
   proposition \ref{Thm.Fast.Loops.via.discriminant.higher.dim}:
   \beq
   p=\frac{p}{q}\cdot q\le (r_k-1)+q\cdot mult(\De_k)\le (r_k-1)+\sum_i q_i\cdot mult(\De_{k,i}).
   \eeq
  It remains to verify: any deformation $\cS(t)$ that preserves the tangent plane, $T_{\cS(t)}=T_{\cS(0)},$ induces
   an equisingular family of planar curve germs, $\De_k(t)\sset \cS(t).$  This is done by explicit computation.
   Denote $mult(\De)=d$ and $mult(\P T_\De, o_k)=m.$ Choose the coordinates in $(\C^n,o)$ such that $ o_k=[0:\dots:0:1]\in \P T_\De,$  but generically otherwise. Then the defining equation of $\De$ belongs to the ideal
      $(x_1,\dots,x_{n-1})^m\cdot (x_1,\dots,x_n)^{d-m}+(x_n)^{d+1}\sset \C\{x\}.$ Moreover, by the genericity of the coordinate choice, the defining equation of $\De$ contains the monomials
    $x^d_1,\dots,x^d_{n-1},$ $x^m_1x^{d-m}_n,\dots,x^m_{n-1}x^{d-m}_n,$ $x^{d+m-1}_n.$  (The last monomial is ensured by the condition $mult(\P T_\De, o_k)> mult(\tDe, o_k).$) After a $GL(n-1)$-transformation on the coordinates $x_1,\dots,x_{n-1},$ the surface $\cS$ is defined by
     the ideal $\bl x_2-a_2(x_1,x_n),\dots,  x_{n-1}-a_{n-1}(x_1,x_n)\br,$ with $a_i\in (x_1)+(x_n)^2.$
      Thus the curve germ $\De\cap \cS$ is defined by : $x^{d+m-k}_n+x^m_1 x^{d-m}_n+x^d_1+(h.o.t.),$ for some  $m>k>0,$ and with some coefficients.
       Any deformation $\cS(t),$ i.e. a deformation of $a_2(x_1,x_n),\dots,   a_{n-1}(x_1,x_n)$ by terms in $(x_1,x_n)^2,$ results in monomials strictly above the Newton diagram. Hence it gives an equisingular family $\De_k(t).$
        Thus both assumptions of proposition \ref{Thm.Fast.Loops.via.discriminant.higher.dim}  are satisfied.
\epr

\bex
Consider hypersurface germ $X=V(x^p_{n+1}-f(x_1,\dots,x_n))\sset (\C^{n+1},o).$ Suppose $ord(f)\ge p\ge 2$ and $f$ is square-free.
 Then the covering $X\to (\C^n_{x_1\dots x_n},o)$ is convenient, and its discriminant, $\De=V(f)\sset (\C^n,o),$ is reduced. The total ramification index over each point of  $\De$ equals $(p-1).$ Therefore the ratio in corollary
  \ref{Thm.Fast.Loops.via.PT.discriminant.higher.dim} is $\frac{p}{q}=\frac{p}{p-1}\le 2.$ Suppose $\De$ is not an ordinary multiple point, i.e.
  $\P T_\De$ has a singular point, say $o_k\in \P^{n-1}.$  If $mult(\P T_\De,o_k)> mult(\tDe,o_k),$ then $X$ has a fast loop.

 As a particular case, let $X=V(x^p_{n+1}-\sum^n_{i=1}x^{d_i}_i)\sset (\C^{n+1},o),$ with $d_1\ge d_2\ge \cdots\ge d_n\ge p\ge2.$
  The discriminant is $\De=V(\sum^n_{i=1}x^{d_i}_i)\sset (\C^n,o).$
  Suppose $d_1>d_n,$ then $\P T_\De$ has multiplicity $d_n$ at the point $o_k:=[1:0:\cdots:0]\in \P^{n-1}\sset Bl_o(\C^n,o).$
   Define the local coordinates at this point via $x_1=\tx_1,x_2=\tx_2\cdot\tx_1,$\dots,  $x_n=\tx_n\cdot\tx_1.$
    Take the strict transform  $\tDe=V(\tx_1^{d_1-d_n}+\cdots+\tx^{d_n}_n)\sset Bl_o(\C^n,o).$
     Then $mult(\P T_\De,o_k)> mult(\tDe,o_k)$ iff $d_1<2d_n.$ Altogether,  if $d_n<d_1<2d_n$ then $X$ has a fast loop.
      In fact, a sufficient condition for a fast loop is just $d_n<d_1,$ see \cite{Kerner-Mendes.Weighted.Homogen}.
\eex

 \end{document}